\documentclass[11pt,twoside,final]{scrartcl}
\usepackage{url,a4wide}
\usepackage{amsmath, amsfonts, amssymb, amsthm, mathtools, nicefrac}
\usepackage{todonotes}
\usepackage{enumitem}
\usepackage{graphicx}

\usepackage{amsmath}
\usepackage{amsthm}
\usepackage{amssymb}
\usepackage{graphics}
\usepackage{showkeys}
\usepackage{pifont,dsfont}
\usepackage{mathrsfs}
\usepackage{euscript}
\DeclareMathAlphabet{\mathpzc}{OT1}{pzc}{m}{it}

\graphicspath{{pics/}}
\usepackage[skip=1pt, font=normalsize]{subcaption}
\usepackage{hyperref} % provides \url command for bibtex and links to jump within documents
\hypersetup{plainpages=false, colorlinks, linkcolor=black, citecolor=black, urlcolor=blue,
pdftitle={Regularized Shannon sampling formulas related to the special affine Fourier transform},
pdfauthor={Frank Filbir, Manfred Tasche, Anna Veselovska},
pdfstartview={FitBH}}
\allowdisplaybreaks

\newcommand{\e}{\mathrm e}

\newcommand{\sinc}{\mathrm{sinc}}

\newcommand{\R}{\mathbb R}
\newcommand{\C}{\mathbb C}
\newcommand{\Z}{\mathbb Z}
\newcommand{\N}{\mathbb N}

\newcommand{\scr}[1]{\ensuremath{\mathscr{#1}}}

\newcommand{\sign}{\mathrm{sign}}
\newcommand{\supp}{\mathrm{supp}}
\newcommand{\im}{\mathrm{i}}
\newcommand{\dx}{\mathrm{d}}

\newcommand{\ds}{\displaystyle}
\newcommand{\ts}{\textstyle}

\renewcommand{\d}[1]{\,\textnormal{d}#1}
\renewcommand{\t}[1]{\,\tilde{#1}}

\newcommand{\banf}{\textsf{Proof.}\ }

\newcommand{\bend}{\hspace*{0ex} \hfill \hbox{\vrule height
    1.5ex\vbox{\hrule width 1.4ex \vskip 1.4ex\hrule  width 1.4ex}\vrule
    height 1.5ex}}
\newtheorem{theorem}{Theorem}[section]
\newtheorem{corollary}[theorem]{Corollary}
\newtheorem{lemma}[theorem]{Lemma}
\newtheorem{alg}[theorem]{Algorithm}
\newtheorem{definition}[theorem]{Definition}
\newtheorem{example}[theorem]{Example}
\newtheorem{remark}[theorem]{Remark}
\newtheorem{proposition}[theorem]{Proposition}

\newenvironment{algorithm}[1]{\goodbreak~\begin{alg}[#1]~\vspace{-9pt}~\\
		\rule{\linewidth}{0.5pt}~\\}{\vspace{-9pt}~\\
		\rule{\linewidth}{0.5pt}~\end{alg}}

\numberwithin{equation}{section}
\numberwithin{table}{section}
\numberwithin{figure}{section}

\renewcommand{\mathbf}[1]{\ensuremath{\boldsymbol{#1}}}

\title{Regularized Shannon sampling formulas related to the special affine Fourier transform}
\author{
	Frank Filbir\footnotemark[1] \and
	Manfred Tasche\footnotemark[3] \and
    Anna Veselovska \footnotemark[2]
}

\date{}

\begin{document}
\maketitle

\begin{abstract}
In this paper, we present new regularized Shannon sampling formulas related to the special affine Fourier transform (SAFT). These sampling formulas use localized
sampling with special compactly supported window functions, namely B-spline, $\sinh$-type, and continuous Kaiser--Bessel window functions. In contrast to the Shannon sampling series for SAFT, the regularized Shannon sampling formulas for SAFT possesses an exponential decay of the
approximation error and are numerically robust in the presence of noise, if certain oversampling condition is fulfilled. Several numerical experiments illustrate the theoretical results.
\medskip
	
	\smallskip
	
	AMS \emph{Subject Classifications}:
	94A20, 42A38, 65T50.
\end{abstract}

\footnotetext[1]{Corresponding author: filbir@helmholtz-muenchen.de, Mathematical Imaging and Data Analysis, Helmholtz Center Munich and Department of Mathematics, Technische Universit\"at M\"unchen, Germany}
\footnotetext[3]{manfred.tasche@uni-rostock.de, Institute of Mathematics, University of Rostock, Germany}
\footnotetext[2]{anna.veselovska@tum.de, Department of Mathematics
\& Munich Data Science Institute, Technical University of Munich, Munich Center for Machine Learning, Garching/Munich, Germany}
\section{Introduction}\label{sec:Intro}
The \emph{special affine Fourier transform} (SAFT) was introduced by S.~Abe and J.T.~Sheridan \cite{AS94} for the study of certain operations on optical
wave functions. For $f \in L^1(\mathbb R)$, it is an integral transform of the form
\begin{equation}\label{eq:Intro1}
\scr F_A f(\omega)=\int_\R f(t)\,  \phi_A(t,\omega)\, \dx t\,, \quad \omega \in \mathbb R\,,
\end{equation}
with the kernel
$$
\phi_A(t,\omega)=\frac{1}{\sqrt{2\pi\, |b|}}\exp\Big[\frac{\mathrm i}{2b}\,\big(a t^2 + 2pt-2\omega t + d \omega^2 + 2(bq -dp)\omega\Big]\,, \quad t,\,\omega \in \mathbb R\,,
$$
depending on a vector $A=(a,b,c,d,p,q)\in\R^6$ which satisfy the conditions $ad-cb=1$ and $b\not= 0$. The name ``special affine Fourier transform'' comes from
the fact that the transform \eqref{eq:Intro1} is related to a special affine transform of the time-frequency coordinates
$$
\begin{pmatrix}
t' \\
\omega'
\end{pmatrix}
=
\begin{pmatrix}
a&b\\
c&d
\end{pmatrix}\,
\begin{pmatrix}
t \\
\omega
\end{pmatrix}
+
\begin{pmatrix}
p\\
q
\end{pmatrix}\,,
$$
where the matrix $\begin{pmatrix}
a&b\\
c&d
\end{pmatrix}$ is an element of the special linear group $SL(2,\R)$, i.e., its determinant is equal to 1. We will not go into a detailed discussion of the
origin and relations of the SAFT to various fields in physics, but we take \eqref{eq:Intro1} merely as a signal transform. The SAFT includes a number of
well-known signal transforms as special cases. Among them are the following:
\begin{itemize}
\item[(1)] For $A=(0,1,-1,0,0,0)$ we get back the classical {\em Fourier transform} $\scr F$, and the set $A=(0,-1,-1,0,0,0)$ gives the {\em inverse Fourier
transform} $\scr F^{-1}$ such that
$$
\scr Ff(\omega)=\frac{1}{\sqrt{2\pi}}\int_\R f(t)\,\e^{-\im\omega t}\d t\,, \quad  \scr F^{-1}f(t)=\frac{1}{\sqrt{2\pi}}\int_\R f(\omega)\,\e^{\im\omega t}\d\omega\,.
$$
For $A=(0,1,-1,0,p,q)\in \R^6$ we obtain the {\em offset Fourier transform}.
\item[(2)] For $A=(a,b,c,d,0,0)\in \R^6$ with $ad-cb=1$ and $abd \not= 0$,  we get the {\em quadratic Fourier transform}
$$
\scr F_Af(\omega)=\frac{1}{\sqrt{2\pi\, |b|}}\int_\R f(t)\, \e^{\frac{\im}{2b}(at^2-2\omega t+d\omega^2)}\, \dx t\,,
$$
which is known in quantum physics as the {\em canonical linear transform}. Special cases of a linear canonical transform are given by
$A=(1,\lambda,0,1,0,0)$ with $\lambda \not= 0$
 ({\em Fresnel transform}), by $A=(\cos(\theta),\sin(\theta),-\sin(\theta),\cos(\theta),0,0)$ with $\theta \notin \pi\,\mathbb Z$ ({\em fractional Fourier transform}), and by $A=(\cosh(\theta),\sinh(\theta),
 \sinh(\theta),\cosh(\theta),0,0)$ with $\theta \not= 0$ ({\em hyperbolic transform}).
For  $A=(\cos(\theta),\sin(\theta),-\sin(\theta),\cos(\theta),p,q) \in \R^6$ with $\theta \notin \pi\,\mathbb Z$, we obtain the {\em offset fractional Fourier transform}.
\end{itemize}
As a matter of fact  the ordinary translation operator $T_x f(t)=f(t-x)$ does not interact nicely with the kernel of the SAFT \eqref{eq:Intro1}, i.e., the function
$T_x \phi_A(t,\omega)=\phi_A(t-x,\omega)$ is in general different from $\phi_A(t,\omega)\, \phi_A(-x,\omega)$ except if $A=(0,1,-1,0,0,0)$. As a consequence, working
with the ordinary translation operator  a number of facts known from ordinary Fourier analysis are no longer valid. This holds in particular for the important Shannon sampling theorem, which allows the reconstruction of a function $f\in L^2(\R)$ with Fourier transform $\scr F f$ supported in $[-\pi,\pi] $
from its samples at integer points using ordinary translates of the cardinal sine function, viz.
\begin{equation}\label{eq:Intro2}
f(t)=\sum_{n\in\Z} f(n)\ \sinc(t-n)\,, \quad t \in \mathbb R\,.
\end{equation}
In case that the SAFT of $f\in L^2(\R)$ is supported on a compact interval, formula \eqref{eq:Intro2} does in general no longer hold. In order to get an analogue
for \eqref{eq:Intro2} in case of the SAFT it is necessary to work with a different concept of translations.
In \cite{FBR22} a generalized translation operator $T^A_x$ was introduced which suits the SAFT and its various consequences for the related harmonic 
analysis were studied. This operator reads as
$$
T^A_x f(t)=\e^{-\im\frac{a}{b}x(t-x)}\, f(t-x)\,, \quad t,\,x \in \mathbb R\,.
$$
It obviously reduces to the ordinary translation in case that $A=(0,1,-1,0,0,0)$.
Moreover, it was demonstrated in \cite{FBR22} that crucial facts from Fourier analysis, as for example the convolution theorem, hold for the SAFT
if we use  with the operator $T^A_x$ instead of $T_x$. In \cite{BZ19} A. Bhandari and A.I. Zayed also studied the SAFT in some detail and worked out a number of aspects of this transform. In \cite{BFR22} the SAFT and related modulation spaces were studied in detail.  Furthermore, in \cite{BZ19} and later in \cite{FBR22} shift invariant spaces related to $T^A_x$ were 
studied and, in particular, an analogue of the Shannon sampling theorem was derived. It states that a function $f\in L^2(\R)$ with 
$\supp(\scr F_Af)\subseteq [-|b|\,\pi+p,|b|\,\pi+p]$ can be reconstructed from its samples on integer points as
\begin{equation}\label{eq:Intro3}
f(t)=\sum_{n\in\Z}\, f(n)\,\e^{-\im\frac{a}{2b}(t^2-n^2)}\,\sinc(t-n)\,, \quad t \in \mathbb R\,.
\end{equation}
Although this representation is quite satisfactory from a theoretical point of view it suffers in the same manner as the ordinary Shannon formula
\eqref{eq:Intro2} from some computational shortcomings.  Apart from the obvious problem of using infinitely many samples, the slow decay of the 
cardinal sine function prevents a good approximation of $f$ by truncation of the series. In order to mitigate these flaws, the cardinal sine function is 
multiplied by a suitable window function $\varphi$. Hence instead of $\sinc$, the function $\varphi\cdot\sinc$ is used for the reconstruction of a function 
$f \in L^2(\mathbb R)$ with $\supp(\scr F_Af)\subseteq [-|b|\,\delta+p,|b|\,\delta+p]$, where $0 < \delta \leq \pi$. In the classical situation, i.e. where 
$A=(0,1,-1,0,0,0)$, localization by a suitable window function and oversampling was used as regularization strategies for the problem. There is a huge stack 
of papers dealing with the problem of a stable and robust computation of the sampling series, see \cite{DDeV03, StTa05, StTa06, PT21, KPT22,KPT23} and references therein.  The Gaussian  
function was frequently used as a window. However, it turned out that certain compactly supported window function lead to much better approximation 
results \cite{KPT22}. Motivated by this observation we concentrate our studies in this paper on a suitable class of compactly supported window functions 
as well.  We consider in particular three types of window functions in detail an analyse there approximation behaviour. 

The paper is organized as follows. In Section 2 we recall all necessary properties of the special affine Fourier transform. Section 3 provides the results on 
sampling for the special affine Fourier transform. Regularization of the Shannon sampling formula will be considered in Section 4. A detailed study of 
regularization by specific window functions will be presented in Section 5. We will complete the paper by Section 6 which is dedicated to numerical 
experiments. 
 
%%%%%%%%%%%%%%%%%%%%%%%%%%%%%%
\section{Special affine Fourier transform}\label{sec:SAFT}
%%%%%%%%%%%%%%%%%%%%%%%%%%%%%%
In this section we shall define the special affine Fourier transform (SAFT) more precisely and  we will discuss the relevant properties of this transform.
Before doing so, let us briefly recall the definition of function spaces which are relevant for us in the sequel. By $C_0(\R)$ we denote the Banach space of 
continuous functions $f:\R\to\C$ which vanish at infinity equipped with the norm $\|f\|_\infty=\max_{t\in\R}|f(t)|$, and $C_c(\R)$ denotes the subspace of 
continuous functions with compact support. The spaces $L^p(\R),\ 1\leq p\leq\infty$ are defined as usual with their respective norm 
$$
\|f\|_p= \begin{cases}
            \Big(\ds\int_\R |f(t)|^p\, \d t\Big)^{1/p}& 1\leq p<\infty,\\
            ess\, \sup_{t\in\R} |f(t)|& p=\infty.
            \end{cases}
$$
The inner product for $L^2(\R)$ is denoted by $\langle\cdot,\cdot\rangle$. \\
The special affine Fourier transform (SAFT) of a function $f\in L^1(\R)$ is defined as 
\begin{equation}\label{eq:SAFT1}
\scr F_A f(\omega) := \frac{1}{\sqrt{2\pi \,|b|}}\, \int_{\mathbb R} f(t)\,
\e^{\frac{\im}{2b}\,\big(a t^2 + 2pt-2\omega t + d \omega^2 + 2(bq -dp)\omega\big)}\, \dx t,\quad \omega\in\R,
\end{equation}
where $A = (a,b,c,d,p,q)\in {\mathbb R}^6$  is a fixed vector with $a d - b c = 1$ and $b\not= 0$. We exclude the limit case $b=0$, as this would make it 
necessary to define the transform in a distributional sense, a setup which we will not consider. For comprehensive studies of the properties of the
SAFT we refer to \cite{BZ19}, \cite{BFR22}, and \cite{FBR22}.\\
 The SAFT can be written in a more convenient form using  the auxiliary functions
\begin{equation}
\label{eq:SAFT1a}
\eta_A(\omega)=\e^{\frac{\im}{2b}(d\omega^2+2(bq-dp)\omega)}\,,\quad \rho_A(t)=\e^{\frac{\im}{2b}(at^2+2pt)}\,, \quad \omega,\,t\in \mathbb R\,.
\end{equation}
The SAFT now reads as
\begin{equation}\label{eq:SAFT2}
\scr F_Af(\omega)=\frac{\eta_A(\omega)}{\sqrt{|b|}}\ \scr F(f\,\rho_A)\big(\frac{\omega}{b}\big)\,, \quad \omega \in \R\,,
\end{equation}
where $\scr F$ stands for the ordinary Fourier transform, given by  \eqref{eq:SAFT1} with $A=(0,1,-1,0,0,0)$. 
As $\eta_A$ and $\rho_A$ are
unimodular functions, i.e., $|\eta_A(\omega)|=1 =|\rho_A(t)|$ for all $t,\,\omega\in\R$, we immediately get from \eqref{eq:SAFT2} that $\scr F_Af$ belongs to
$C_0(\R)$ with $\|\scr F_Af\|_\infty\leq (2\pi |b|)^{-1/2}\,\|f\|_1$. Moreover, \eqref{eq:SAFT2} also shows that $\scr F_A$ can be extended to $L^2(\R)$ and
defines a unitary operator on that space, viz.
\begin{equation}\label{eq:SAFT2a}
\langle\scr F_Af,\scr F_A g\rangle=\langle f, g\rangle,\quad f,\,g\in L^2(\R).
\end{equation}

In particular, $\|\scr F_Af\|_2=\|f\|_2$ for $f\in L^2(\R)$. The inverse of $\scr F_A$ on $L^2(\R)$ can readily be obtained from \eqref{eq:SAFT2}. It is given as 
\begin{equation}\label{eq:SAFT3}
\scr F_A^{-1} f(t)=\frac{\bar{\rho}_A(t)}{\sqrt{2\pi|b|}}\int_\R f(\omega) \, \bar{\eta}_A(\omega)\, \e^{\im\omega t/b}\, \dx \omega.
\end{equation}
In \cite{FBR22}, the authors introduced a generalized translation operator related to the SAFT which was called $A$-{\em translation}. 
For a function $f:\R\to\C$ and $x\in\R$ it reads as
\begin{equation}\label{eq:SAFT4}
T^A_xf(t)=\e^{-\im\frac{a}{b} x(t-x)}\, f(t-x),\quad t\in\R.
\end{equation}
Obviously, $T^A_x$ is norm preserving on all spaces $L^p(\R),\ 1\leq p\leq \infty$, and it reduces to the ordinary translation $T_xf(t)=f(t-x)$ if 
$A=(0,1,-1,0,0,0)$. However, note that $T^A_x(f\cdot g)$ is in general different from $(T^A_xf)\cdot(T^A_x g)$. The following statements regarding the 
generalized translation operator $T^A_x$  can be found in \cite{FBR22}. The proofs of these statements are rather straightforward and we will omit them 
here.
\begin{proposition}\label{pro:SAFT1}
For the $A$-translation $T^A_x$ and any $f\in L^1(\R)$, the following properties hold
\begin{itemize}
\item[(i)] $T^A_x\, T^A_y=\e^{-\im\frac{a}{b} x y}\, T^A_{x+y},\quad x,\,y\in\R\,,$
\item[(ii)] $\scr F_A(T^A_x f)(\omega)=\rho_A(x)\,\e^{-\im\omega x/b}\, \scr F_Af(\omega),\quad x,\,\omega\in\R\,.$
\end{itemize}
\end{proposition}
With help of the $A$-translation we define the $A$-{\em convolution} of two functions $f,\,g\in L^1(\R)$ as
$$
(f\star_A g)(t)=\frac{1}{\sqrt{2\pi|b|}}\int_\R f(x)\, T^A_xg(t)\ \dx x\,.
$$
Then we have  $f\star_A g\in L^1(\R)$ and  $\| f\star_A g\|_1\leq (2\pi |b|)^{-1/2} \|f\|_1 \|g\|_1$. Moreover, 
$$
\scr F_A(f\star_A g)(\omega)=\bar{\eta}_A(\omega)\, \scr F_Af(\omega)\, \scr F_A g(\omega)\,.
$$
Again, the latter two equations reduce to the classical expression in case $A=(0,1,-1,0,0,0)$, i.e.,
$$
(f\star g)(t)=\frac{1}{\sqrt{2\pi}}\int_\R f(x)\, g(t-x)\ \dx x\,, \quad \scr F(f\star g)(\omega)=\scr F f(\omega)\, \scr F g(\omega).
$$
The proofs of these facts can be found in \cite{FBR22}.  \\
It is convenient to define to introduce the {\em chirp modulation operator}
\begin{equation}
\label{eq:SAFT4a}
C_sf(t)=\e^{\frac{\im}{2} st^2} f(t)\,, \quad s,\,t \in\R\,.
\end{equation}
The A-translation operator and the chirp modulation operator do obviously not commute. We have 
\begin{equation}\label{eq:SAFT4b}
T_x^A C_sf(t)=\e^{-\frac{\im}{2}(t-x)[(2\frac{a}{b}+s)x-s\,t]} f(t-x),
\end{equation}
and if $s=-\frac{a}{b}$ we get 
\begin{equation}\label{eq:SAFT5}
T^A_x\, C_{-\frac{a}{b}}\,f(t)=\e^{-\im\frac{a}{2b} (t^2-x^2)}\, f(t-x)\,.
\end{equation}

\section{Shannon sampling series for SAFT}\label{sec:SSS}
In \cite{FBR22}, an analogue of the Shannon sampling theorem for the SAFT was obtained. In order to formulate the sampling theorem, let us introduce some
notation. For $0<\delta\leq\pi$ let $I_\delta=[-\delta|b|,\,\delta|b|]$ and define the subspace $B^A_{p,\delta}$ of $L^2(\R)$ as
$$
B^A_{p,\delta}=\{f\in L^2(\R): \supp(\scr F_A f)\subseteq p+I_\delta\}\,.
$$
This space is a generalization of the space of ordinary $\delta$-bandlimited functions $B_\delta=\{f\in L^2(\R):\supp(\scr F f)\subseteq [-\delta,\,\delta]\}$ to the case
of the SAFT. Between the spaces $B_\delta$ and $B_{p,\delta}^A$, the following close relations hold:

\begin{proposition}\label{pro:SSS1}
If $g \in B_\delta$, then $f = C_{-\frac{a}{b}}g$ belongs to $B^A_{p,\delta}$.\\
If $f \in B_{p,\delta}^A$, then it holds $g = C_{\frac{a}{b}}f \in B_\delta$.
\end{proposition}

\banf Assume that $g\in B_\delta$. Using \eqref{eq:SAFT1a}, \eqref{eq:SAFT2}, and \eqref{eq:SAFT4a}, we obtain
\begin{eqnarray}
\label{eq:SSS1a}
\scr F_A(C_{-\frac{a}{b}}g)(\omega)&=&\frac{\eta_A(\omega)}{\sqrt{|b|}}\,\scr F\big((C_{-\frac{a}{b}}g)\, \rho_A\big)\big(\frac{\omega}{b}\big) \nonumber \\
&=& \frac{\eta_A(\omega)}{\sqrt{2\pi\,|b|}}\,\int_{\mathbb R} g(t)\,{\mathrm e}^{- {\mathrm i}\,\frac{\omega - p}{b}\,t}\,{\mathrm d}t \nonumber \\
&=& \frac{\eta_A(\omega)}{\sqrt{|b|}}\,(\scr F g)(\textstyle{\frac{\omega-p}{b}})\,.
\end{eqnarray}
From $g \in B_\delta$ it follows that $(\scr F g)(\textstyle{\frac{\omega-p}{b}}) = 0$ for all $\omega \in \mathbb R$ with $|\omega - p| > \delta\,|b|$, i.e.,
for all $\omega \in \mathbb R \setminus (p + I_\delta)$. Hence for $f = C_{-\frac{a}{b}}g$ equation \eqref{eq:SSS1a} gives 
$$
\mathrm{supp}\,(\scr F_A f) \subseteq p + I_\delta\,,
$$
that means $f \in B^A_{p,\delta}$.\\
Now we consider $f \in B^A_{p,\delta}$. Applying \eqref{eq:SAFT1a} and \eqref{eq:SAFT4a}, we conclude
\begin{eqnarray*}
\scr F(f \,\rho_A)\big(\frac{\omega}{b}\big) &=& \frac{1}{\sqrt{2\pi}}\,\int_{\mathbb R} f(t)\,\rho_A(t)\,{\mathrm e}^{- {\mathrm i}\,\frac{\omega}{b}\,t}\d t \\
&=&  \frac{1}{\sqrt{2\pi}}\,\int_{\mathbb R} (C_{\frac{a}{b}}f)(t)\,{\mathrm e}^{- {\mathrm i}\,\frac{\omega-p}{b}\,t}\,{\mathrm d}t \\
&=& \scr F(C_{\frac{a}{b}}f)\big(\frac{\omega-p}{b}\big)\,.
\end{eqnarray*}
Then from \eqref{eq:SAFT2} it follows that
\begin{equation}\label{eq:SSS1c}
(\scr F_A f)(\omega)= \frac{\eta_A(\omega)}{\sqrt{|b|}}\,\scr F(f\,\rho_A)\big(\frac{\omega}{b}\big)=\frac{\eta_A(\omega)}{\sqrt{|b|}}\,\scr F(C_{\frac{a}{b}}f)\big(\frac{\omega-p}{b}\big)\,.
\end{equation}
By assumption it holds
$$
\mathrm{supp}\,(\scr F_A f) \subseteq p + I_\delta\,.
$$
Hence by \eqref{eq:SSS1c} and $|\eta_A (\omega)| = 1$ for all $\omega \in \mathbb R$, we obtain
$$
\mathrm{supp}\,\scr F(C_{\frac{a}{b}}f)\big(\frac{\cdot -p}{b}\big)\subseteq p + I_\delta\,,
$$
that means
$$
\mathrm{supp}\,\scr F(C_{\frac{a}{b}}f) \subseteq [-\delta,\, \delta]\,.
$$
Thus $g = C_{\frac{a}{b}}f$ belongs to $B_\delta$. \\
Therefore the mapping $\scr I:\,B_\delta \to B_{p,\delta}^A$ with $\scr I g = C_{-\frac{a}{b}} g$ for each $g \in B_\delta$
is a linear isometry, since we have
$$
\| \scr I g\|_2 = \|C_{-\frac{a}{b}} g\|_2 = \| g \|_2
$$
for all $g \in B_\delta$. 
\bend
\medskip

Note that the equivalence class of functions $f\in L^2(\R)$ with $\supp(\scr F(C_{-\frac{a}{b}}f))\subseteq [-\delta,\delta]$ always contains a smooth function.
Indeed, for any $r\in\N_0$ the function $(\im\omega)^r\, \scr F(C_{-\frac{a}{b}}f)(\omega)$ is in $L^1([-\delta,\delta])$ which implies that
$$
(C_{-\frac{a}{b}}f)^{(r)}(t)=\frac{1}{\sqrt{2\pi}}\int_{-\delta}^\delta \scr F(C_{-\frac{a}{b}}f)(\omega)\, (\im\omega)^r\,\e^{\im\omega t}\d\omega
$$
belongs to $C_0(\R)$. Clearly, with $C_{-\frac{a}{b}}f$ also $f$ is smooth. In the sequel we will always pick the smooth representative from the equivalence 
class $C_{-\frac{a}{b}}f$. \\

If $f \in  B^A_{p,\delta}$ with $0 < \delta < \pi$, we will speak about \emph{oversampling} of $f$ by $\{f(n)\}_{n\in \Z}$.
The cardinal sine function is defined as usual by
\begin{equation}
\label{eq:SSS1b}
\sinc(t)= \begin{cases}
  \frac{\sin(\pi t)}{\pi t}\,,  & t\in\R\setminus\{0\}\,, \\
   1\,,& t=0\,.
\end{cases}
\end{equation}
Define $\psi(t)=C_{-\frac{a}{b}} \sinc\,(t)$. As $\sinc\in B_\pi$ we have $\psi\in B^A_{p,\delta}$ in view of Proposition \ref{pro:SSS1}.  The following 
Shannon sampling theorem for the SAFT was proved in \cite{BZ19} and in \cite{FBR22}. 
For some special SAFT's, corresponding sampling theorems 
were studied mainly in the signal processing literature. In particular sampling theorems for the Fresnel transform \cite{Go81}, fractional Fourier transform 
\cite{ZG99, BM10, BZ12}, and the linear canonical transform \cite{ZRMT09, TLWA08,SLSZ12} were studied before. However, they are all straightforward 
consequences of following theorem.
\begin{theorem}\label{theo:SSS1}
If $f\in B^A_{p,\delta}$, then it holds
\begin{equation}\label{eq:SSS1}
f(t)=\sum_{n\in\Z} f(n)\, T^A_n\, \psi(t)
\end{equation}
for every $t\in\R$.
\end{theorem}
Using the relation \eqref{eq:SAFT5}, the sampling series \eqref{eq:SSS1} can be rewritten as
\begin{equation}\label{eq:SSS2}
f(t)=\sum_{n\in\Z} f(n)\, \e^{-\im\frac{a}{2b}(t^2-n^2)}\, \sinc(t-n).
\end{equation}
In order to obtain the convergence  of the truncated sampling series
$$
S_Nf(t)=\sum_{n=-N}^N f(n)\, T^A_n\psi(t)\,, \quad t \in \R\,,
$$
to the function $f\in B^A_{p,\delta}$ for $N \to \infty$, we first prove the following lemma.
\begin{lemma}\label{lem:SSS1}
Let $\delta\in(0,\pi]$ and $f\in B^A_{p,\delta}$ be given. Then we have
\begin{equation}\label{eq:SSS2a}
\|f\|_2^2=\sum_{n\in\Z} |f(n)|^2.
\end{equation}
\end{lemma}
\banf Since
$$
\langle T^A_n\psi, T^A_m\psi\rangle=\langle T_n\,\sinc, T_m\,\sinc\rangle=\delta_{n,m}
$$
for all $n,\,m\in\Z$, it follows that $\{T^A_n\psi:n\in\Z\}$ is an orthonormal system in $L^2(\R)$. If $f\in B^A_{p,\delta}$ we get from Theorem \ref{theo:SSS1}
$$
f(t)=\sum_{n\in\Z} f(n)\, T^A_n\psi(t)\,.
$$
This leads to
$$
\langle f,f\rangle=\Big\langle\sum_{n \in \Z} f(n)\, T^A_n\psi, \sum_{m \in \Z} f(m)\, T^A_m\psi\Big\rangle=\sum_{n\in \Z} |f(n)|^2.
$$
\bend\\
We are now prepared to prove the desired convergence result for the truncated sampling series.

\begin{proposition}\label{pro:SSS2}
Let $\delta\in (0,\pi]$ and $f\in B^A_{p,\delta}$ be given. Then we have
\begin{equation}\label{eq:SSS2b}
\lim_{N\to\infty}\, \|f-S_Nf\|_{\infty}=0\,.
\end{equation}
Moreover, for each $t\in\R$ the sampling series \eqref{eq:SSS1} converges absolutely and we have
\begin{equation}\label{eq:SSS3}
\sum_{n\in\Z} |f(n)|\, |\sinc(t-n)|\leq \|f\|_2\,.
\end{equation}
\end{proposition}
\banf Using  Theorem \ref{theo:SSS1} we get
$$
f(t)-S_Nf(t)=\sum_{|n|>N} f(n)\, T^A_n\psi(t)\,.
$$
The Cauchy-Schwarz inequality yields
\begin{eqnarray}\label{eq:SSS4}
\big|f(t) - S_N f(t) \big| &\le& \Big(\sum_{|n|>N} |f(n)|^2\Big)^{1/2}\, \Big(\sum_{|n|>N} \big| T_n^A \psi(t)\big|^2\Big)^{1/2} \nonumber \\
&\le& \Big(\sum_{|n|>N} |f(n)|^2\Big)^{1/2}\, \Big(\sum_{|n|>N} \big|\sinc(t - n)\big|^2 \Big)^{1/2}\,.
\end{eqnarray}
Since
\begin{equation}\label{eq:SSS5}
\sum_{n\in\Z}|\sinc(t-n)|^2=1\,, \quad t \in \mathbb R\,,
\end{equation}
the second factor of \eqref{eq:SSS4} converges to zero  uniformly with respect to $t\in\R$ as $N \to \infty$. In view of
\eqref{eq:SSS2a} the first factor converges to zero for $N \to \infty$. This implies \eqref{eq:SSS2b}.\\
The statement about the absolute convergence of the sampling series \eqref{eq:SSS1} resp. \eqref{eq:SSS2} and the estimate \eqref{eq:SSS3}  follow again
by applying the Cauchy-Schwarz inequality, \eqref{eq:SSS2a}, and \eqref{eq:SSS5}.
\bend

Although, the reconstruction of $f\in B^A_{p,\delta}$ from samples $f(n), n\in\Z$, is possible according to Theorem \ref{theo:SSS1}, it is not a stable process.
In order to demonstrate this consider $f\in B^A_{p,\delta}$ with $\delta\in(0,\pi]$. In view of Proposition \ref{pro:SSS1} we may without lost of generality assume that $f$
is smooth. For a sufficiently large $N\in\N$ define erroneous samples of $f(n)$, $n\in\Z$, by setting
$$
\tilde{f}(n)=
\begin{cases}
f(n)+\varepsilon_n& -N\leq n\leq N\,,\\
f(n) & \mbox{otherwise}\,,
\end{cases}
$$
with error terms $\varepsilon_n\in\C$ and $|\varepsilon_n|\leq \varepsilon$. We obtain
$$
\begin{array}{ll}
\tilde{f}(t)&=\ds\sum_{n\in\Z} \tilde{f}(n)\, \e^{-\im\frac{a}{2b}(t^2-n^2)}\, \sinc(t-n)\\[3ex]
&=f(t)+\ds\sum_{n=-N}^N\varepsilon_n\, \e^{-\im\frac{a}{2b}(t^2-n^2)}\, \sinc(t-n)\,,\quad t\in\R\,.
\end{array}
$$
Thus an upper estimate is given by
$$
\|\tilde{f}-f\|_{\infty}\leq \varepsilon\,\max_{t\in\R}\sum_{n=-N}^N|\sin(t-n)|
$$
By \cite[Theorem 2.2]{KPT22} we have
$$
\max_{t\in\R}\sum_{n=-N}^N|\sinc(t-n)|<\frac{2}{\pi}[\,\ln(N)+2\ln(2)+\gamma\,]+\frac{N+2}{\pi N(N+1)}\,,
$$
with Euler's constant \cite{Y91}
$$
\gamma:=\lim_{N\to\infty}\Big(\sum_{n=1}^N\frac{1}{n}-\ln(N)\Big)=0.57721566\dots\, .
$$
We shall now give a lower bound for the approximation of $f$ utilizing perturbed sample values. The result shows in particular that the reconstruction of
$f\in B^A_{p,\delta}$ by the Shannon sampling series \eqref{eq:SSS1} resp. \eqref{eq:SSS2} is numerically unstable.

\begin{theorem}\label{theo:SSS2}
Let $\delta\in(0,\pi]$ and $f\in B^A_{p,\delta}$ be given. For $N\in\N$ and $\varepsilon>0$ define
$$
\varepsilon_n=\varepsilon\ \sign\big(\sinc({\ts\frac{1}{2}}-n)\big)\,\e^{-\im\frac{a}{8b}(4n^2-1)}=\varepsilon\,(-1)^{n+1}\sign(2n-1)\,\e^{-\im\frac{a}{8b}(4n^2-1)}
$$
for $-N\leq n\leq N$ and $\varepsilon_n=0$ for $|n|>N$. Then it holds
$$
\|\tilde{f}-f\|_{\infty}\geq \varepsilon\,\big(\frac{2}{\pi}\ln(N)+\frac{4}{\pi}\,\ln(2)+\frac{2\gamma}{\pi}\big)>\varepsilon\big(\frac{2}{\pi}\ln(N)+\frac{5}{4}\big)\,.
$$
\end{theorem}
\banf
The special choice of the error term leads to
$$
\tilde{f}(t)-f(t)=\varepsilon\,\sum_{n=-N}^N\sign\big(\sinc({\ts\frac{1}{2}}-n)\big)\,\e^{-\im\frac{a}{8b}(4n^2-1)}\,\e^{-\im\frac{a}{2b}(n^2-t^2)}\,\sinc(t-n)\,.
$$
For $t=\frac{1}{2}$ we obtain
$$
\|\tilde{f}-f\|_{\infty}\geq |\tilde{f}({\ts\frac{1}{2})-f(\frac{1}{2})}|=\varepsilon\, \sum_{k=-N}^N |\sinc(t-n)|=\frac{2\varepsilon}{(2N+1)\pi}+\frac{4\varepsilon}{\pi}
\sum_{n=1}^N\frac{1}{2n-1}\,.
$$
From \cite[Theorem 2.2]{KPT22} we get
$$
\sum_{n=1}^N\frac{1}{2n-1}> \frac{1}{2}\ln(N)+\ln(2)+\frac{\gamma}{2}-\frac{1}{4N(2N+1)},
$$
which yields
$$
\|\tilde{f}-f\|_{\infty}> \varepsilon\,\big(\frac{2}{\pi}\ln(N)+\frac{4}{\pi}\ln(2)+\frac{2\gamma}{\pi}\big)+\varepsilon\,\frac{2N+1}{\pi\, N(2N+1)}
>\varepsilon\,\big(\frac{2}{\pi}\ln(N)+\frac{4}{\pi}\ln(2)+\frac{2\gamma}{\pi}\big).
$$
Since $\frac{4}{\pi}\ln(2)+\frac{2\gamma}{\pi}=1.2500093\dots>\frac{5}{4}$, we get the final estimate.
\bend
\section{Regularized Shannon sampling formulas for SAFT}\label{sec:RSS}
The practical use of the Shannon sampling series \eqref{eq:SSS2} for SAFT is rather limited due to the need of infinitely many samples.  A second limitation of its use in
practice originates in the fact that the cardinal sine function decays very slowly, which results in poor convergence of the sampling series. Furthermore, in case
of noisy samples $f(n)$ the convergence of the Shannon sampling series can even break down completely (see \cite{DDeV03, KPT23}). In order to overcome these limitations, regularization
techniques were considered. For the classical Fourier transform, the concept of regularized Shannon sampling formulas with localized sampling and oversampling has been studied
by several authors \cite{Q03, Q04, Q06, SchSt07, TaSuMu07, LZ16, KPT22, KPT23}. Often a Gaussian window function supported on whole $\mathbb R$ was used in these studies.
In \cite{KPT22}, it was shown that the compactly supported $\sinh$-type window function \eqref{eq:RSS3} produces smaller approximation errors than the Gaussian window function.
Therefore we focus on compactly supported window functions in the following.

For any $m\in\N\setminus\{1\}$ let $\Phi_m$ be the class of even continuous functions $\varphi:\R\to [0,1]$ with the following properties:
\begin{itemize}
\item[(i)] $\varphi$ is supported on $[-m,m]$,
\item[(ii)] $\varphi$ is monotonically decreasing on $[0,m]$ with $\varphi(0)=1$,
\item[(iii)] the Fourier transform of $\varphi$ is explicitly known.
\end{itemize}
A function from the class $\Phi_m$ will henceforth be called {\em window function}. The following window functions will be considered in the remaining part 
of the paper.
\begin{example}\label{ex:RSS1}
\begin{itemize}
\item[(a)] For $m\in\N$ denote by $M_{m}$ the centered cardinal $\mathrm B$-spline of order $m$ which is defined recursively by 
$$
M_m(t)=\int_{-\frac{1}{2}}^{\frac{1}{2}} M_{m-1}(t-\tau)\d\tau,\ m\geq1,\quad M_1(t)=\mathds{1}_{[-\frac{1}{2},\frac{1}{2}]}.
$$
Note that $M_m$ is supported on $[-\frac{m}{2},\frac{m}{2}]$. For even  order $2s\in 2\N$ we define the $\mathrm B$-{\em spline window function} as
\begin{equation}\label{eq:RSS3}
\varphi_{\mathrm B}(t)=\frac{1}{M_{2s}(0)}\,M_{2s}\Big(\frac{s\,t}{2\pi m}\Big)\,.
\end{equation}
\item[(b)] The $\sinh$-{\em type window function} is defined for a positive parameter $\beta$ as
\begin{equation}\label{eq:RSS1}
\varphi_{\sinh}(t)=
\begin{cases}
 \frac{1}{\sinh\beta}\sinh\big(\beta\sqrt{1-(\frac{t}{2\pi m})^2}\big)\,,  & t\in [-2\pi\, m, 2\pi\, m]\,,\\
 0\,,& \mbox{elsewhere}\,.
\end{cases}
\end{equation}
\item[(c)] Let again $\beta$ be a positive number. The {\em continuous Kaiser-Bessel window function} is defined as
\begin{equation}\label{eq:RSS2}
\varphi_{\mathrm{cKB}}(t)=
\begin{cases}
 \frac{1}{I_0(\beta)-1}\Big( I_0\big(\beta\sqrt{1-(\frac{t}{2\pi m})^2}\big)-1\Big)\,,  & t\in [-2\pi\, m,2\pi\, m]\,,\\
 0\,,& \mbox{elsewhere}\,,
\end{cases}
\end{equation}
where $I_0$ is the modified Bessel function of the first kind, given by \cite{Wa44}
$$
I_0(x)=\sum_{k=0}^\infty \frac{1}{((2k)!!)^2}\, x^{2k}\,, \quad x \in \R\,.
$$
\end{itemize}
\end{example}

\begin{figure}[ht]
    \centering 
    \subfloat[\label{fig:-1}]{
        \includegraphics[width=0.45\linewidth]{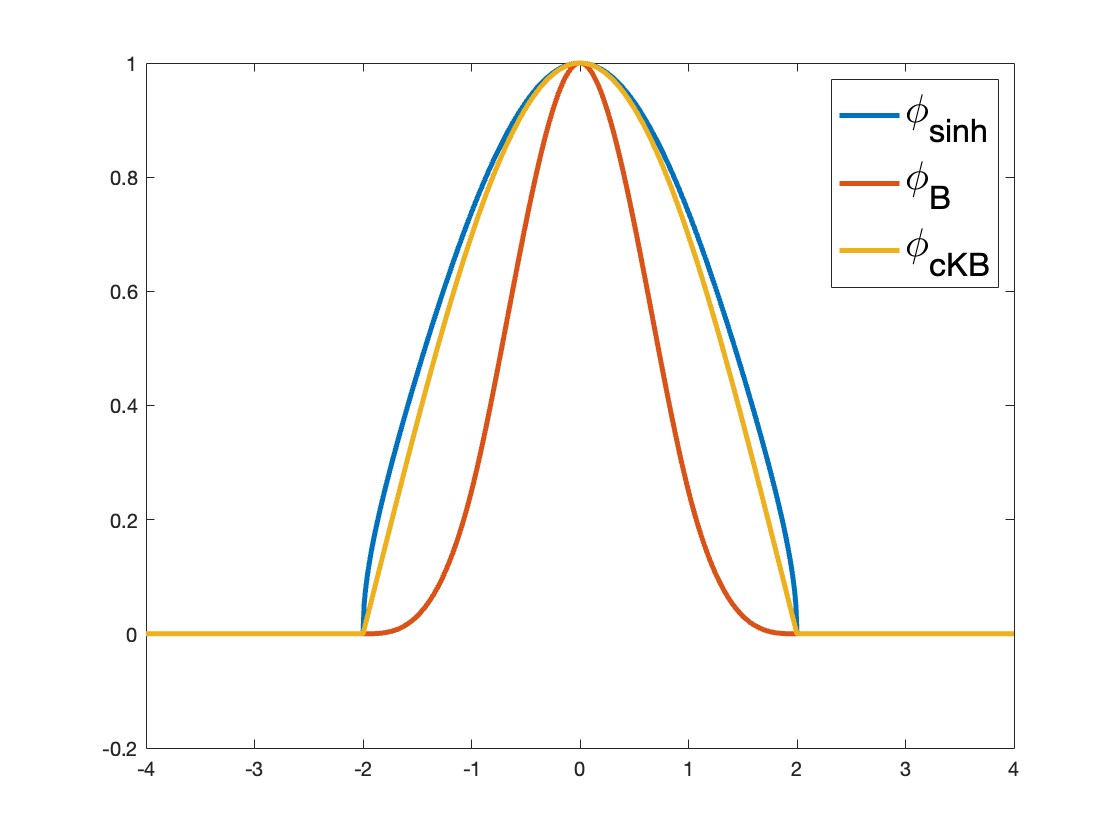}
    }
    \subfloat[\label{fig:-2}]{
       \includegraphics[width=0.45\linewidth]{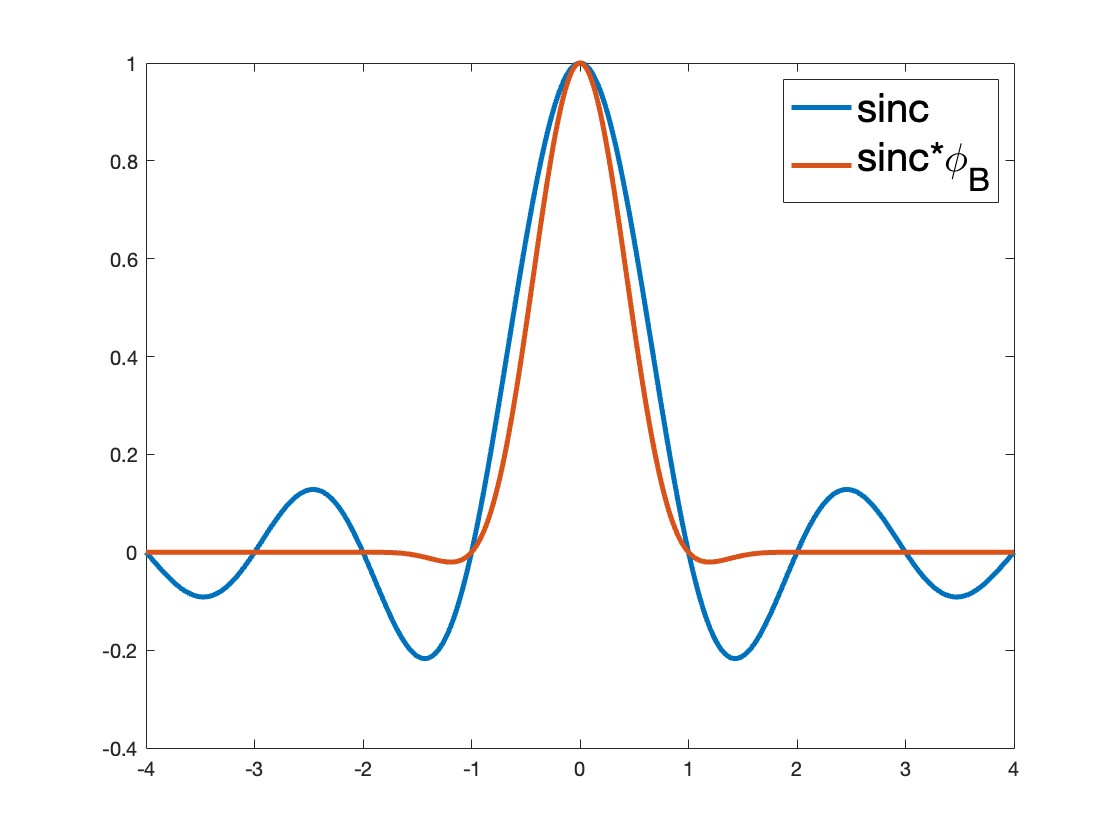}
    }
    \caption[]{(a) Window functions from Example~\ref{ex:RSS1} for the window size $m=2$.  \\(b) Regularization of the sinc function via the centered cardinal B-spline window.}
    \label{fig:1:window-fun}
\end{figure}

Particular cases of the above-indicated window functions and a regularized version of the sinc function are depicted in Figure~\ref{fig:1:window-fun}. 

We will consider the following regularization strategies for Shannon sampling series of~SAFT:
\begin{itemize}
\item[1.] The sampling series \eqref{eq:SSS2} will be regularized by a window function $\varphi \in \Phi_m$, viz., we consider instead of \eqref{eq:SSS2} the {\em regularized Shannon sampling formula for the} SAFT $\scr F_A$ given by
\begin{equation}\label{eq:RSS0}
R^A_{\varphi,m} f(t)=\sum_{n\in\Z} f(n)\,\e^{-\im\frac{a}{2b}(t^2-n^2)}\,\sinc(t-n)\, \varphi(t-n)\,, \quad t \in \R\,.
\end{equation}
Since $\varphi$ is compactly supported, the computation of $R^A_{\varphi,m} f(t)$ for fixed $t\in \R$ needs only a finite number of samples of $f$.
\item[2.] Another established technique is to oversample the function $f$, for which we need the stronger condition $f\in B^A_{p,\delta}$ with some $\delta\in (0,\, \pi)$.
\end{itemize}
Note that since 
$$
{\mathrm e}^{- {\mathrm i}\,\frac{a}{2b}\,(t^2 - n^2)}\,\mathrm{sinc}(t-n)\,\varphi(t-n)\big|_{t=k} = \delta_{n,k}\,, \quad n,\,k\in \mathbb Z\,.
$$
we have the following interpolation property 
$$
R^A_{\varphi,m}f(k)=f(k)\,, \quad k \in \mathbb Z\,,
$$
for $f\in B_{p,\delta}^A$. \\
Furthermore, note that the use of a window function $\varphi\in\Phi_m$ implies that the computation of $R^A_{\varphi,m}f(t)$ for
$t\in\R\setminus\Z$ requires only $2m$ samples $f(n)$, where $n\in\Z$ fulfills $|n-t|< m$. Hence the function $f$ can be recovered on the interval $[0,\,1]$ by
$$
f(t) =
\begin{cases}
 f(0)\,,  & t = 0\,,\\
 f(1)\,,  & t = 1\,,\\
\ds \sum_{n=1-m}^m f(n)\,\e^{-\im\frac{a}{2b}(t^2-n^2)}\,\sinc(t-n)\, \varphi(t-n)\,, & t \in (0,\,1)\,.
\end{cases}
$$
Thus the reconstruction of $f$ on the interval $[-N,N]$
with $N\in\Z$ needs only $2N+2m-1$ samples $f(n)$ with $n=1-N-m,\,\dots,\, N+m-1$.

Now we give an estimate for the uniform approximation error $\| f-R^A_{\varphi,m}f\|_{\infty}$ for the approximation of $f$ by $R^A_{\varphi,m}f$.

\begin{theorem}\label{theo:RSS1}
Let $f\in B^A_{p,\delta}$ with $\delta\in (0,\pi]$ be given. Further let $\varphi\in\Phi_m$. Then we have
$$
\| f-R^A_{\varphi,m}f\|_{\infty}\leq E(m,\delta)\ \|f\|_2\,,
$$
with
\begin{equation}\label{eq:RSS4}
E(m,\delta)=\sqrt{\frac{\delta}{\pi}}\ \max_{\omega\in[-\delta,\delta]}\Big| 1-\frac{1}{\sqrt{2\pi}}\int_{\omega-\pi}^{\omega+\pi}\hat{\varphi}(\tau)\, \dx\tau\Big|.
\end{equation}
\end{theorem}
\banf
We shall compute the SAFT of the error function
$$
e(t)=f(t)-R^A_{\varphi,m}f(t)\quad t\in\R.
$$
Applying $\scr F_A$ to the representation \eqref{eq:SSS1} and using  Proposition \ref{pro:SAFT1} (ii)  leads to
$$
\scr F_A f(\omega)=\sum_{n\in\Z} f(n)\, \rho_A(n)\,\e^{-\im n\omega/b}\, \scr F_A\psi(\omega)
$$
with $\psi(t) = C_{-\frac{a}{b}}\,\mathrm{sinc}(t)$.
A straightforward computation shows
$$
\scr F_A\psi(\omega)=\frac{1}{\sqrt{2\pi|b|}}\,\eta_A(\omega)\, \mathds{1}_{[-\pi,\pi]}\big(\frac{\omega-p}{b}\big)
$$
with the indicator function $\mathds{1}_{[-\pi,\pi]}$ of the interval $[-\pi,\,\pi]$. Hence it holds
$$
\scr F_A f(\omega)=\frac{1}{\sqrt{2\pi|b|}}\,\mathds{1}_{[-\pi,\pi]}\big(\frac{\omega-p}{b}\big)\, \sum_{n\in\Z} f(n)\, \rho_A(n)\,\e^{-\im n\omega/b}.
$$
For the computation of the SAFT of $R^A_{\varphi,m}f$, let
$$
\chi(t)= C_{-\frac{a}{b}}\big(\varphi\, \sinc\big)(t)=\e^{-\im\frac{a}{2b}t^2}\,\varphi(t)\, \sinc(t) \,.
$$
The SAFT of the function $\chi$ can be computed using \eqref{eq:SSS1a},
we obtain
$$
\scr F_A\chi(\omega)=\frac{1}{\sqrt{2\pi\,|b|}}\,\eta_A(\omega) \int_\R \sinc(t)\, \varphi(t)\, \e^{-\im(\omega-p)t/b}\, \dx t=\frac{1}{\sqrt{|b|}}\,\eta_A(\omega)\
\scr F(\varphi\, \sinc)\big(\frac{\omega-p}{b}\big).
$$
The convolution property of $\scr F$ in $L^2(\R)$ (see \cite[Lemma 2.2]{KPT22}) gives  $(\scr F\varphi)\ast(\scr F\sinc) = \scr F(\varphi\cdot\sinc)$, 
and since $\scr F\sinc=\frac{1}{\sqrt{2\pi}}\,\mathds{1}_{[-\pi,\pi]}$ and $\scr F\varphi=\hat{\varphi}$, we get
$$
\scr F(\varphi\cdot\sinc)(\omega)=\frac{1}{\sqrt{2\pi}}\,\big(\hat{\varphi}\ast \mathds{1}_{[-\pi,\pi]}\big)(\omega)
= \frac{1}{2\pi}\,\int_{\omega-\pi}^{\omega+\pi}\hat{\varphi}(\tau)\,\dx\tau\,.
$$
Hence we arrive at
$$
\scr F_A(R^A_{\varphi,m}f)(\omega)=\frac{\eta_A(\omega)}{2\pi\sqrt{|b|}}\int_{-\pi+(\omega-p)/b}^{\pi+(\omega-p)/b}\hat{\varphi}(\tau)\, \dx \tau\
\sum_{n\in\Z} f(n)\,\rho_A(n)\,\e^{-\im n\omega/b}\,.
$$
Thus, with the auxiliary function
\begin{equation}\label{eq:RSS4a}
\Delta(\omega)=\mathds{1}_{[-\pi,\pi]}(\omega)-\frac{1}{\sqrt{2\pi}}\int_{-\pi+\omega}^{\pi+\omega}\hat{\varphi}(\tau)\,\dx\tau\,, \quad \omega \in \R\,,
\end{equation}
the SAFT of the error function $e$ can be written as
$$
\scr F_Ae(\omega)=\scr F_A f(\omega)\, \Delta\big(\frac{\omega-p}{b}\big)\,.
$$
Taking its inverse SAFT we finally obtain the following estimate
$$
\begin{array}{ll}
|e(t)|&=\ds\Big|\frac{1}{\sqrt{2\pi|b|}}\bar{\rho}_A(t)\ds\int_\R\scr F_A e(\omega)\, \bar{\eta}_A(\omega)\, \e^{\im\omega t/b}\, \dx\omega\Big|\\[4ex]
&\leq\ds \frac{1}{\sqrt{2\pi|b|}}\int_\R|\scr F_A f(\omega)|\, \big|\Delta\big(\frac{\omega-p}{b}\big)\big|\, \dx\omega.\\[4ex]
\end{array}
$$
Taking into account that $f\in B^A_{p,\delta}$ we arrive at
\begin{equation}\label{eq:RSS5}
|e(t)|\leq \frac{1}{\sqrt{2\pi|b|}}\int_{p-|b|\,\delta}^{p+|b|\, \delta}\scr |\scr F_A f(\omega)|\,\dx\omega\
\max_{\omega'\in[-\delta,\delta]}\Big| 1-\frac{1}{\sqrt{2\pi}}\int_{-\pi+\omega'}^{\pi+\omega'} \hat{\varphi}(\tau)\, \dx\tau\Big|
\end{equation}
with $\omega'=(\omega-p)/b\in[-\delta,\delta]$. The Cauchy-Schwarz inequality in combination with the Parseval equality \eqref{eq:SAFT2a} gives
$$
\int_{p-|b|\,\delta}^{p+|b|\, \delta}\scr |\scr F_A f(\omega)|\,\dx\omega\leq \sqrt{2\delta\,|b|}\, \|\scr F_A f\|_2=\sqrt{2\delta\,|b|}\, \|f\|_2\,.
$$
With this  inequality \eqref{eq:RSS5} finally yields
$$
\|e\|_{\infty}\leq \sqrt{\frac{\delta}{\pi}}\ \|f\|_2\, \max_{\omega'\in[-\delta,\delta]}\Big| 1-\frac{1}{\sqrt{2\pi}}
\int_{-\pi+\omega'}^{\pi+\omega'} \hat{\varphi}(\tau)\, \dx\tau\Big|\,.
$$
This completes the proof. \bend
\medskip

In the previous section we have demonstrated that the Shannon sampling series does not behave stably with respect to perturbation of the samples of the
function. Now we are considering this problem for the regularized Shannon sampling formula related to SAFT. It turns out that in contrast to the sampling series \eqref{eq:SSS1} resp.
\eqref{eq:SSS2} the regularized Shannon sampling formula \eqref{eq:RSS1} for SAFT is numerically robust, i.e., the uniform error
$\|R^A_{\varphi,m}\tilde{f}-R^A_{\varphi,m}f\|_{\infty}$ is small for small perturbations of the samples of the function $f$. More precisely we have the
following statement.
\begin{theorem}\label{theo:RSA2}
For $\delta\in(0,\pi]$ let $f\in B^A_{p,\delta}$, and let $\varphi\in\Phi_m$ be a window function. Furthermore, let $\tilde{f}(n)=f(n)+\varepsilon_n,\, n\in\Z$, be noisy samples,
where the complex sequence $\{\varepsilon_n\}_{n\in\Z}$ satisfies $|\varepsilon_n|\leq\varepsilon$ for some $\varepsilon>0$.\\
Then the following estimates
hold
\begin{eqnarray}
%\begin{array}{rl}
\|R^A_{\varphi,m}\tilde{f}- R^A_{\varphi,m}f\|_{\infty}&\leq &\varepsilon\, \big(2+\sqrt{2\pi}\,\hat{\varphi}(0)\big)\,,\label{eq:RSS6}\\[2ex]
\|f-R^A_{\varphi,m}\tilde{f}\|_{\infty}&\leq &\|f- R^A_{\varphi,m}f\|_{\infty}+\varepsilon\, \big(2+\sqrt{2\pi}\,\hat{\varphi}(0)\big)\,.\label{eq:RSS7}
%\end{array}
\end{eqnarray}
\end{theorem}
\banf We consider the perturbation error
$$
\t e(t) =R^A_{\varphi,m}\t f(t)-R^A_{\varphi,m}f(t)=\sum_{n\in\Z} \varepsilon_n\,\e^{\im\frac{a}{2b}(n^2-t^2)}\,\sinc(t-n)\,\varphi(t-n).
$$
Note that for any $t\in\R$ the sum contains only finitely many non-vanishing terms. We first consider the case $t\in(0,\,1)$. Due to the properties of
$\varphi\in\Phi_m$ and the fact $|\varepsilon_n| \leq \varepsilon$, we obtain
$$
|\t e(t)|\leq \sum_{n=1-m}^m |\varepsilon_n|\,|\sinc(t-n)|\,\varphi(t-n)\leq \varepsilon\,\sum_{n=1-m}^m \varphi(t-n)\,.
$$
Since $\varphi\in\Phi_m$ is monotonically decreasing on $[0,m]$ we have
$$
\begin{array}{ll}
\ds\sum_{n=1-m}^m \varphi(t-n)&=\ds\Big(\sum_{n=1-m}^0+\sum_{n=1}^m\Big)\varphi(t-n)=\sum_{n=0}^{m-1}\varphi(t+n)+\sum_{n=1}^m\varphi(t-n)\\[4ex]
&\leq \ds\sum_{n=0}^{m-1}\varphi(n)+\sum_{n=1}^m\varphi(1-n)=2\sum_{n=0}^{m-1}\varphi(n)\,.
\end{array}
$$
The latter sum can be estimated further by applying once more the monotonicity of $\varphi$
$$
\sum_{n=0}^{m-1}\varphi(n)<\varphi(0)+\int_0^{m-1}\varphi(t)\d t\leq\varphi(0)+\int_0^m\varphi(t)\d t\,.
$$
Since $\hat{\varphi}(0)=\sqrt{\frac{2}{\pi}}\int_0^m\varphi(t)\d t$ we eventually have
$$
|\t e(t)|\leq2\,\varepsilon\sum_{n=0}^{m-1}\varphi(n)\leq\varepsilon\,\big(2\,\varepsilon\,\varphi(0)+\sqrt{2\pi}\,\hat{\varphi}(0)\big)
=\varepsilon\,\big(2+\sqrt{2\pi}\,\hat{\varphi}(0)\big)
$$
for every $t\in(0,\,1)$. Note that $|\t e(0)|=|\varepsilon_0|\leq \varepsilon$ and $|\t e(1)|=|\varepsilon_1|\leq\varepsilon$, which leads to
\begin{equation}\label{eq:RSS8}
\max_{t\in[0,1]}|\t e(t)|\leq \varepsilon\,\big(2+\sqrt{2\pi}\,\hat{\varphi}(0)\big)\,.
\end{equation}
The same technique can obviously be applied to get the estimate \eqref{eq:RSS8} for every interval $[k,k+1],\, k\in\Z$, and thus \eqref{eq:RSS7} follows. \\
The inequality \eqref{eq:RSS7} can easily be derived by applying the triangular inequality and \eqref{eq:RSS6}.
\bend
\section{Specific regularized Shannon sampling formulas for SAFT}\label{sec:SRSS}
In this section we will study regularizations of the Shannon sampling series by specific window functions which were introduced in Example \ref{ex:RSS1}.

\subsection{Regularization with B-spline window function}\label{ssec:BSR}
First we consider the B-spline window function as defined in \eqref{eq:RSS3} with $m\in\N\setminus\{1\}$ and $s =\lceil\frac{m+1}{2}\rceil$.
Assume that a function $f \in B^A_{p,\delta}$ with $0 <\delta < \pi - \frac{2s}{m}$ is given. Then the
B-{\em spline regularized Shannon sampling formula for} $\scr F_A$
{\em and} $f$ reads as follows
\begin{equation}\label{eq:SRSS1}
R^A_{{\mathrm B},m}f(t)=\sum_{n\in\Z}f(n)\,\e^{-\im\frac{a}{2b}(t^2-n^2)}\,\sinc(t-n)\,\varphi_{\mathrm B}(t-n)\,.
\end{equation}
We will now show that the approximation error $\|f-R^A_{{\mathrm B},m}f\|_{\infty}$ decays exponentially with respect to $m$. More precisely the following statement
holds.
\begin{theorem}\label{theo:SRSS1}
Let $\varphi_{\mathrm B}\in\Phi_m$ be the ${\mathrm B}$-spline window function \eqref{eq:RSS3} with $m \in \N\setminus \{1\}$ and
$s=\lceil\frac{m+1}{2}\rceil$.
Let $f\in B^A_{p,\delta}$ with $0 < \delta < \pi - \frac{2s}{m}$ be given.\\
Then it holds
\begin{equation}\label{eq:SRSS2}
\|f-R^A_{{\mathrm B},m}f\|_{\infty}\leq \frac{1}{\pi}\,\Big(\frac{2s}{m\,(\pi-\delta)}\Big)^{2s-1}\, \|f\|_2\,.
\end{equation}
\end{theorem}
\banf According to Theorem \ref{theo:RSS1} we have
$$
\|f-R^A_{{\mathrm B},m}f\|_{\infty}\leq\sqrt{\frac{\delta}{\pi}}\ \|f\|_2\, \max_{\omega\in[-\delta,\delta]}
\Big|1-\frac{1}{\sqrt{2\pi}}\int_{\omega-\pi}^{\omega+\pi} \hat{\varphi}_{\mathrm B}(\tau)\d\tau\Big|\,.
$$
By \cite[p. 452]{PPST18} we have
$$
\int_\R M_{2s}(t)\, \e^{-\im\omega t}\d t=\big(\sinc\,\frac{\omega}{2}\big)^{2s}\,,
$$
since the cardinal sine function is defined by \eqref{eq:SSS1b}.
Hence the B-spline window function \eqref{eq:RSS3} has the Fourier transform
\begin{equation}
\label{eq:hatvarphiB}
\hat{\varphi}_{\mathrm B}(\omega)=\frac{m}{\sqrt{2\pi}\,s\,M_{2s}(0)}\,\big(\sinc\,\frac{m\,\omega}{2s}\big)^{2s}\,,
\end{equation}
which results in
$$
1=\varphi_{\mathrm B}(0)=\frac{1}{\sqrt{2\pi}}\int_\R\hat{\varphi}_{\mathrm B}(\tau)\d \tau=\frac{m}{2\pi\,s\,M_{2s}(0)}\int_\R\big(\sinc\,\frac{m\,\tau}{2s}\big)^{2s}\d\tau\,.
$$
Then for $\omega\in[-\delta,\delta]$ the auxiliary function \eqref{eq:RSS4a} takes the form
$$
\begin{array}{ll}
\Delta(\omega)&=\ds\frac{m}{2\pi\,s\,M_{2s}(0)}\Big[\int_\R\big(\sinc\,\frac{m\,\tau}{2s}\big)^{2s}\d\tau-
\int_{\omega-\pi}^{\omega+\pi}\big(\sinc\,\frac{m\,\tau}{2s}\big)^{2s}\d\tau\Big]\\[4ex]
&=\ds\frac{m}{2\pi\,s\,M_{2s}(0)}\Big[\int_{\pi-\omega}^\infty\big(\sinc\,\frac{m\,\tau}{2s}\big)^{2s}\d\tau+
\int_{\omega+\pi}^{\infty}\big(\sinc\,\frac{m\,\tau}{2s}\big)^{2s}\d\tau\Big].
\end{array}
$$
For the single integral terms we have the following estimates
$$
\int_{\pi\pm\omega}^\infty\big(\sinc\,\frac{m\,\tau}{2s}\big)^{2s}\d\tau\leq\frac{(2s)^{2s}}{m^{2s}}\int_{\pi\pm\omega}^\infty\tau^{-2s}\d\tau
=\frac{(2s)^{2s}}{(2s-1)\,m^{2s}(\pi\pm\omega)^{2s-1}}
$$
which gives
$$
\Delta(\omega)\leq \frac{(2s)^{2s-1}}{(2s-1)\,m^{2s-1}\,\pi\,M_{2s}(0)}\Big[\frac{1}{(\pi-\omega)^{2s-1}}+\frac{1}{(\pi+\omega)^{2s-1}}\Big]
$$
for all $\omega\in[-\delta,\delta]$.  Since $0<\delta<\pi- \frac{2s}{m}$ and $\omega\in[-\delta,\delta]$ we have $\pi-\omega,\,\pi+\omega\in[\pi-\delta,\,\pi+\delta]$. Taking into
account that the function $x^{1-2s}$ is decreasing for $x>0$, we conclude
$$
\max_{\omega\in [-\delta,\,\delta]}|\Delta(\omega)|\leq \frac{2\,(2s)^{2s}}{(2s-1)\,m^{2s-1}\,\pi\,M_{2s}(0)\,(\pi-\delta)^{2s-1}}\,.
$$
Applying \cite[formula (5.3)]{KPT22} gives 
\begin{equation}\label{eq:SRSS3}
\frac{4}{3}\leq\sqrt{2s}\, M_{2s}(0)<\sqrt{\frac{6}{\pi}}\,,
\end{equation}
and thus we obtain for the error term \eqref{eq:RSS4} the estimate
$$
E(m,\delta)\leq\frac{3}{2\pi}\sqrt{\frac{\delta}{\pi}}\ \Big(\frac{2s}{m(\pi-\delta)}\Big)^{2s-1}\,\frac{\sqrt{2s}}{2s-1} \leq \frac{1}{\pi}\sqrt{\frac{\delta}{\pi}}\ \Big(\frac{2s}{m(\pi-\delta)}\Big)^{2s-1}\,.
$$
The assumption $s=\lceil\frac{m+1}{2}\rceil\geq 2$ implies $2s-2\leq m\leq 2s-1$. By the oversampling condition  $0 < \delta < \pi - \frac{2s}{m}$ it holds $0 < \frac{2s}{m \,(\pi -\delta)}<1$.
This completes the proof. \bend
\medskip

Now we show that for the B-spline regularized Shannon sampling formula \eqref{eq:SRSS1} for $\scr F_A$ and $f \in B_{p,\delta}^A$, the uniform perturbation error
$\| R_{{\mathrm B},m}^A f - R_{{\mathrm B},m}^A \tilde f\|_{\infty}$ only grows as ${\mathcal O}(\sqrt m)$, if the oversampling condition $0 < \delta < \pi - \frac{2s}{m}$ is fulfilled.

\begin{theorem}\label{theo:Brobust}
Let $\varphi_{\mathrm B}\in\Phi_m$ be the ${\mathrm B}$-spline window function \eqref{eq:RSS3} with $m \in \N\setminus \{1\}$ and
$s=\lceil\frac{m+1}{2}\rceil$.
Let $f\in B^A_{p,\delta}$ with $0 < \delta < \pi - \frac{2s}{m}$ be given. Furthermore, let $\tilde{f}(n)=f(n)+\varepsilon_n,\, n\in\Z$, be noisy samples,
where the complex sequence $\{\varepsilon_n\}_{n\in\Z}$ satisfies $|\varepsilon_n|\leq\varepsilon$ for some $\varepsilon>0$.\\
Then the $\mathrm B$-spline regularized Shannon sampling formula \eqref{eq:SRSS1} for $\scr F_A$ and $f \in B_{p,\delta}^A$ is numerically robust and it holds
$$
\|R^A_{{\mathrm B},m}f - R_{{\mathrm B},m}^A \tilde f\|_{\infty}\leq \varepsilon\, \big(2 + \frac{3}{2}\, \sqrt m\big)\,.
$$
\end{theorem}

\banf By Theorem \ref{theo:RSA2} we only have to determine the value $\hat \varphi_{\mathrm B}(0)$ for the B-spline window function \eqref{eq:RSS3}. From \eqref{eq:hatvarphiB} and \eqref{eq:SRSS3} it follows that
$$
\hat \varphi_{\mathrm B}(0) = \frac{m}{\sqrt{2 \pi}\,s \,M_{2s}(0)} \leq \frac{3\,m}{4\,\sqrt{s\,\pi}}\,.
$$
Due to $s =\lceil\frac{m+1}{2}\rceil$ it holds $\sqrt s \geq \sqrt{\frac{m}{2}}$ and hence
$$
\sqrt{2 \pi}\,\hat \varphi_{\mathrm B}(0)  \leq \frac{3}{2}\,\sqrt m\,.
$$
By the oversampling condition $0 < \delta < \pi- \frac{2s}{m}$, the error $\|f - R_{{\mathrm B},m}^A \tilde f\|_{\infty}$ is small too by Theorem \ref{theo:SRSS1} and the triangle inequality. \bend

\subsection{Regularization with sinh-type window function}\label{ssec:STR}
Now we consider the regularization of the Shannon sampling series for the SAFT $\scr F_A$ and  $f\in B^A_{p,\delta}$ with $\delta\in(0,\pi)$ using the $\sinh$-type window function
\eqref{eq:RSS1} with parameter $\beta=m(\pi-\delta)$. The $\sinh$-{\em type regularized Shannon sampling formula for} $\scr F_A$ {\em and} $f$ takes the form
\begin{equation}\label{eq:SRSS4}
R^A_{\sinh,m}f(t)=\sum_{n\in\Z} f(n)\,\e^{-\im\frac{a}{2b}(t^2-n^2)}\,\sinc(t-n)\, \varphi_{\sinh}(t-n)\,, \quad t \in \R\,.
\end{equation}
We will demonstrate that the uniform error of approximating $f$ by $R^A_{\sinh,m}f$ decays exponentially with respect to $m$.
\begin{theorem}\label{theo:SRRS2}
Let $\delta\in (0,\pi)$ and $f\in B^A_{p,\delta}$ be given. Further let $\varphi_{\sinh}$ be the $\sinh$-type window function as defined in \eqref{eq:RSS1} with parameter
$\beta=m(\pi-\delta)$. Then it holds
$$
\| f-R^A_{\sinh,m}f\|_{\infty}\leq \sqrt{\frac{\delta}{\pi}}\,\e^{-m(\pi-\delta)}\, \|f\|_2\,.
$$
\end{theorem}

\banf According to Theorem \ref{theo:RSS1} we have
$$
\| f-R^A_{\sinh,m}f\|_{\infty}\leq \sqrt{\frac{\delta}{\pi}}\,\| f \|_2\,  \max_{\omega\in[-\delta,\delta]}|\Delta(\omega)|
$$
with
$$
\Delta(\omega)=1-\frac{1}{\sqrt{2\pi}}\int_{\omega-\pi}^{\omega+\pi}\hat{\varphi}_{\sinh}(\tau)\d\tau\,, \quad \omega \in [-\delta,\,\delta]\,.
$$
Following \cite[p.~38, 7.58]{Ob90} we have
\begin{equation}\label{eq:SRSS5}
\hat{\varphi}_{\sinh}(\tau)=\frac{m\sqrt{\pi}}{\sqrt{2}\sinh\beta}\cdot
\begin{cases}
(1-\nu^2)^{-1/2}\, I_1(\beta\sqrt{1-\nu^2})\,,& |\nu|<1\,,\\[2ex]
(\nu^2-1)^{-1/2}\, J_1(\beta\sqrt{\nu^2-1})\,,& |\nu|>1
\end{cases}
\end{equation}
with $\nu=\frac{m}{\beta}\tau$.  With this change of variables the function $\Delta$ now reads as
$$
\Delta(\omega)=1-\frac{\beta}{\sqrt{2\pi}\,m}\int_{-\nu_1(-\omega)}^{\nu_1(\omega)} \hat{\varphi}_{\sinh}\big(\frac{\beta}{m}\nu\big)\d\nu
$$
with the linear increasing function $\nu_1(\omega)=\frac{m}{\beta}(\omega+\pi),\ \omega\in[-\delta,\delta]$. Since $\delta\in(0,\pi)$ and $\beta=m(\pi-\delta)$,
we have $\nu_1(-\delta)=1$ and $\nu_1(\omega)\geq 1$ for all $\omega\in[-\delta,\delta]$.\\
Now let
$\Delta(\omega)=\Delta_1(\omega) -\Delta_2(\omega)$ with
$$
\begin{array}{ll}
\Delta_1(\omega)&=1-\ds\frac{\beta}{\sqrt{2\pi}\,m}\int_{-1}^1 \hat{\varphi}_{\sinh}\big(\frac{\beta}{m}\nu\big)\d\nu\,,\\[4ex]
\Delta_2(\omega)&=\ds\frac{\beta}{\sqrt{2\pi}\,m}\Big(\int_{-\nu_1(-\omega)}^{-1}+\int_1^{\nu_1(\omega)}\Big)
   \hat{\varphi}_{\sinh}\big(\frac{\beta}{m}\nu\big)\d\nu\\[4ex]
&=\ds\frac{\beta}{\sqrt{2\pi}\,m}\Big(\int_1^{\nu_1(-\omega)}+\int_1^{\nu_1(\omega)}\Big) \hat{\varphi}_{\sinh}\big(\frac{\beta}{m}\nu\big)\d\nu\,.
\end{array}
$$
In view of \eqref{eq:SRSS5} these functions take the form
$$
\begin{array}{ll}
\Delta_1(\omega)&=1-\ds\frac{\beta}{2\sinh\beta}\int_{-1}^1 \frac{I_1(\beta\sqrt{1-\nu^2})}{\sqrt{1-\nu^2}}\d\nu\,,\\[4ex]
\Delta_2(\omega)&=\ds\frac{\beta}{2\sinh\beta}\Big(\int_1^{\nu_1(-\omega)}+\int_1^{\nu_1(\omega)}\Big) \frac{J_1(\beta\sqrt{\nu^2-1})}{\sqrt{\nu^2-1}}\d\nu\,.
\end{array}
$$
Using \cite[6.681--3]{GR80} and \cite[10.2.13]{abst}, we get
$$
\int_{-1}^1 \frac{I_1(\beta\sqrt{1-\nu^2})}{\sqrt{1-\nu^2}}\d\nu=\int_{-\pi/2}^{\pi/2} I_1(\beta\cos\,s)\d s=\pi\Big(I_{1/2}\big(\ts\frac{\beta}{2}\big)\Big)^2=\frac{4}
{\beta}\Big(\sinh \frac{\beta}{2}\Big)^2
$$
and hence
$$
\Delta_1(\omega)=1-\frac{2\,(\sinh\frac{\beta}{2})^2}{\sinh\beta}=\frac{2\,\e^{-\beta}}{1+\e^{-\beta}}.
$$
By \cite[6.645--1]{GR80} we have
$$
\int_1^\infty\frac{J_1(\beta\sqrt{\nu^2-1})}{\sqrt{\nu^2-1}}\d\nu=I_{1/2}\big({\ts\frac{\beta}{2}}\big)\,K_{1/2}\big({\ts\frac{\beta}{2}}\big)
=\frac{1-\e^{-\beta}}{\e^{\beta}}\,,
$$
where $I_{1/2}$ and $K_{1/2}$ are modified Bessel functions of half order (see \cite[10.2.13, 10.2.14, and 10.2.17]{abst}). Numerical experiments have shown that 
$$
0 < \int_1^T \frac{J_1\big(\beta\,\sqrt{\nu^2-1}\big)}{\sqrt{\nu^2-1}}\,{\mathrm d}\nu \le \frac{3\,\big(1 - {\mathrm e}^{-\beta}\big)}{2\,\beta}\,.
$$
for all $T > 1$.
Thus we obtain
$$
0 \le \Delta_2(\omega) \le \frac{\beta}{2\,\sinh \beta}\,\frac{3\,\big(1 - {\mathrm e}^{-\beta}\big)}{\beta} = \frac{3\,{\mathrm e}^{-\beta}}{1 + {\mathrm e}^{-\beta}}\,.
$$
Altogether the estimate
$$
|\Delta(\omega)\big| \le \big|  \Delta_1(\omega) - \Delta_2(\omega) \big| \le \frac{{\mathrm e}^{-\beta}}{1 + {\mathrm e}^{-\beta}}< {\mathrm e}^{-\beta}
$$
holds for all $\omega \in [-\delta,\, \delta]$. This implies 
$$
E(m,\delta) \le \sqrt{\frac{\delta}{\pi}}\,{\mathrm e}^{-m\,(\pi-\delta)},
$$
which completes the proof.
\bend

Now we show that for the $\sinh$-type regularized Shannon sampling formula \eqref{eq:SRSS4} with respect to $\scr F_A$ and $f \in B_{p,\delta}^A$, the uniform perturbation
error $\|R_{\sinh,m}^A f - R_{\sinh,m}^A \tilde f\|_{\infty}$ only grows as ${\mathcal O}(m)$, if the oversampling condition $0 < \delta \leq \pi - \frac{\pi}{m}$ is fulfilled.

\begin{theorem}\label{theo:sinhrobust}
For $m \in \N\setminus \{1\}$ and $\delta \in \big(0,\,\pi - \frac{\pi}{m}\big]$, let $f \in B_{p,\delta}^A$ be given.
Let $\varphi_{\sinh}\in\Phi_m$ be the $\sinh$-type window function \eqref{eq:RSS1} with $\beta = m\,(\pi - \delta)$.
Furthermore, let $\tilde{f}(n)=f(n)+\varepsilon_n,\, n\in\Z$, be noisy samples,
where the complex sequence $\{\varepsilon_n\}_{n\in\Z}$ satisfies $|\varepsilon_n|\leq\varepsilon$ for some $\varepsilon>0$.\\
Then the $\sinh$-type regularized Shannon sampling formula \eqref{eq:SRSS4} for $\scr F_A$ and $f \in B_{p,\delta}^A$ is numerically robust and it holds
$$
\|R^A_{\sinh,m}f - R_{\sinh,m}^A \tilde f\|_{\infty}\leq \varepsilon\, \big(2 + \frac{3}{2}\, m\big)\,.
$$
\end{theorem}

\banf By Theorem \ref{theo:RSA2} we only have to determine the value $\hat \varphi_{\sinh}(0)$ for the $\sinh$-type window function \eqref{eq:RSS1}. From \eqref{eq:SRSS5} it follows that
$$
{\hat \varphi}_{\sinh}(0) = \frac{m\,\sqrt \pi\,I_1(\beta)}{\sqrt 2  \,\sinh \beta}\,.
$$
Applying the inequality $\sqrt{2 \pi \beta}\,{\mathrm e}^{-\beta}\,I_1(\beta) < 1$ (see \cite[Lemma 7]{PT21}), we find that
$$
{\hat \varphi}_{\sinh}(0) < \frac{m\,{\mathrm e}^{\beta}}{2\,\sqrt \beta  \,\sinh \beta} = \frac{\sqrt m}{\sqrt{\pi - \delta}\,\big(1 - {\mathrm e}^{-2 \beta}\big)}\,.
$$
Under the oversampling condition $0 < \delta \leq \pi - \frac{\pi}{m}$ with $m\geq 2$ it holds $\beta = m\,(\pi-\delta)\geq \pi$ and hence
$$
\frac{1}{\sqrt{\pi - \delta}} \leq \frac{\sqrt m}{\sqrt \pi}\,.
$$
Therefore by Theorem \ref{theo:RSA2} we can estimate
$$
\|R^A_{\sinh,m}f - R_{\sinh,m}^A \tilde f\|_{\infty}\leq \varepsilon\, \big(2 + \frac{\sqrt 2}{1 - {\mathrm e}^{-2 \beta}}\, m\big)\,.
$$
From $\beta \geq \pi$ it follows that
$$
\frac{\sqrt 2}{1 - {\mathrm e}^{-2 \beta}} \leq \frac{\sqrt 2}{1 - {\mathrm e}^{-2 \pi}} < \frac{3}{2}\,.
$$
This completes the proof. \bend

\subsection{Regularization with continuous Kaiser-Bessel window function}\label{ssec:cKB}
Finally we consider the regularization of the Shannon sampling series for the SAFT $\scr F_A$ and  $f\in B^A_{p,\delta}$ with $\delta\in(0,\pi)$ using the continuous Kaiser-Bessel window function
\eqref{eq:RSS2} with parameter $\beta=m(\pi-\delta)$. The {\em continuous Kaiser-Bessel regularized Shannon sampling formula for} 
$\scr F_A$ {\em and} $f$ takes the form
\begin{equation}\label{eq:cKB}
R^A_{\mathrm{cKB},m}f(t)=\sum_{n\in\Z} f(n)\,\e^{-\im\frac{a}{2b}(t^2-n^2)}\,\sinc(t-n)\, \varphi_{\mathrm{cKB}}(t-n)\,, \quad t \in \R\,.
\end{equation}
Now we show that the uniform approximation error $\|f - R^A_{\mathrm{cKB},m}f\|_{\infty}$ decays exponentially with respect to $m$.
\begin{theorem}\label{theo:cKB}
Let $\delta\in (0,\pi)$ and $f\in B^A_{p,\delta}$ be given. Further let $\varphi_{\mathrm{cKB}}$ be the continuous Kaiser-Bessel window function 
\eqref{eq:RSS2} with parameter $\beta=m(\pi-\delta)$. Then 
$$
\| f-R^A_{\mathrm{cKB},m}f\|_{\infty}\leq \sqrt{\frac{\delta}{\pi}}\,\frac{1}{I_0(\beta)-1}\, \|f\|_2.
$$
\end{theorem}

\banf From Theorem \ref{theo:RSS1} it follows that
$$
\| f-R^A_{\mathrm{cKB},m}f\|_{\infty}\leq \sqrt{\frac{\delta}{\pi}}\,\| f \|_2\,  \max_{\omega\in[-\delta,\delta]}|\Delta(\omega)|
$$
with
$$
\Delta(\omega)=1-\frac{1}{\sqrt{2\pi}}\int_{\omega-\pi}^{\omega+\pi}\hat{\varphi}_{\mathrm{cKB}}(\tau)\d\tau\,, \quad \omega \in [-\delta,\,\delta]\,.
$$
According to  \cite[p.~3, 1.1, and p.~95, 18.31]{Ob90}, the Fourier transform of \eqref{eq:RSS2} has the form
\begin{equation}\label{eq:hatcKB}
\hat{\varphi}_{\mathrm{cKB}}(\tau)=\frac{\sqrt{2}\,m}{\sqrt{\pi}\,(I_0(\beta)-1)}\cdot
\begin{cases}
\Big( \frac{\sinh(\beta\,\sqrt{1 - \nu^2})}{\beta\, \sqrt{1 - \nu^2}} - \frac{\sin(\beta \nu)}{\beta \nu}\Big)\,,& 0<|\nu|<1\,,\\[4ex]
\Big( \frac{\sin(\beta\,\sqrt{\nu^2-1})}{\beta\, \sqrt{\nu^2-1}} - \frac{\sin(\beta \nu)}{\beta \nu}\Big)\,,& |\nu|>1,
\end{cases}
\end{equation}
with $\nu=\frac{m}{\beta}\tau$. With this substitution the function $\Delta$ now reads as
$$
\Delta(\omega)=1-\frac{\beta}{\sqrt{2\pi}\,m}\int_{-\nu_1(-\omega)}^{\nu_1(\omega)} \hat{\varphi}_{\mathrm{cKB}}\big(\frac{\beta}{m}\nu\big)\d\nu
$$
with the linear increasing function $\nu_1(\omega)=\frac{m}{\beta}(\omega+\pi)$ for $\omega\in[-\delta,\delta]$. Since $\delta\in(0,\pi)$ and $\beta=m(\pi-\delta)$,
we have $\nu_1(-\delta)=1$ and $\nu_1(\omega)\geq 1$ for all $\omega\in[-\delta,\delta]$.\\
Now we split
$\Delta(\omega)=\Delta_1(\omega) -\Delta_2(\omega)$ with
$$
\begin{array}{ll}
\Delta_1(\omega)&=1-\ds\frac{\beta}{\sqrt{2\pi}\,m}\int_{-1}^1 \hat{\varphi}_{\mathrm{cKB}}\big(\frac{\beta}{m}\nu\big)\d\nu\,,\\[4ex]
\Delta_2(\omega)&=\ds\frac{\beta}{\sqrt{2\pi}\,m}\Big(\int_{-\nu_1(-\omega)}^{-1}+\int_1^{\nu_1(\omega)}\Big)
   \hat{\varphi}_{\mathrm{cKB}}\big(\frac{\beta}{m}\nu\big)\d\nu\\[4ex]
&=\ds\frac{\beta}{\sqrt{2\pi}\,m}\Big(\int_1^{\nu_1(-\omega)}+\int_1^{\nu_1(\omega)}\Big) \hat{\varphi}_{\mathrm{cKB}}\big(\frac{\beta}{m}\nu\big)\d\nu\,.
\end{array}
$$
Using \eqref{eq:hatcKB}, these functions take the form
$$
\begin{array}{ll}
\Delta_1(\omega)&=1-\ds\frac{\beta}{\pi\,(I_0(\beta)-1)}\int_{-1}^1 \Big( \frac{\sinh(\beta\,\sqrt{1 - \nu^2})}{\beta\, \sqrt{1 - \nu^2}} - \frac{\sin(\beta \nu)}{\beta \nu}\Big) \d\nu\,,\\[4ex]
\Delta_2(\omega)&=\ds\frac{\beta}{\pi \,(I_0(\beta)-1)}\Big(\int_1^{\nu_1(-\omega)}+\int_1^{\nu_1(\omega)}\Big) \Big( \frac{\sin(\beta\,\sqrt{\nu^2-1})}{\beta\, \sqrt{\nu^2-1}} - \frac{\sin(\beta \nu)}{\beta \nu}\Big)\d\nu\,.
\end{array}
$$
By \cite[3.997--1]{GR80} we have
\begin{align*}
	\int_{-1}^{1} \frac{\sinh\big(\beta \sqrt{1 - \nu^2}\,\big)}{\beta \sqrt{1 - \nu^2}}\,{\mathrm d}\nu &= \frac{2}{\beta}\, \int_0^1 \frac{\sinh\big(\beta \sqrt{1 - \nu^2}\,\big)}{\sqrt{1 - \nu^2}}\,{\mathrm d}\nu \\[2ex]
	&= \frac{2}{\beta}\, \int_0^{\pi/2} \sinh(\beta \cos \sigma)\,{\mathrm d}\sigma = \frac{\pi}{\beta}\,{\textbf L}_0(\beta)\,,
\end{align*}
where \mbox{${\textbf L}_0$} denotes the \emph{modified Struve function} given by (see \cite[12.2.1]{abst})
$$
	{\textbf L}_0(x) = \sum_{k=0}^{\infty} \frac{(x/2)^{2k+1}}{\big(\Gamma\big(k + \tfrac{3}{2}\big)\big)^2} = \frac{2 x}{\pi}\,\sum_{k=0}^{\infty}\frac{x^{2k}}{\big((2k+1)!!\big)^2}\,, \quad x \in \R\,.
$$
Note that the function \mbox{$I_0(x) - {\textbf L}_0(x)$} is completely monotonic on \mbox{$[0,\,\infty)$} (see \cite[Theorem~1]{BP14}) and tends to zero as \mbox{$x\to \infty$}.
Applying the \emph{sine integral function}
$$
	\mathrm{Si}(x) = \int_0^{x} \frac{\sin t}{t}\,{\mathrm d}t \,, \quad x \in \mathbb R\,,
$$
implies
$$
	\int_{-1}^1 \frac{\sin(\beta \nu)}{\beta \nu}\,{\mathrm d}\nu = 2\,\int_0^1 \frac{\sin(\beta \nu)}{\beta \nu}\,{\mathrm d}\nu  = \frac{2}{\beta}\,\mathrm{Si}(\beta)\,.
$$
Hence we obtain
\begin{align*}
	\Delta_1(\omega) &= 1 - \frac{1}{I_0(\beta) - 1}\,\bigg({\textbf L}_0(\beta) - \frac{2}{\pi}\,\mathrm{Si}(\beta)\bigg) %\\
	= \frac{1}{I_0(\beta) - 1}\,\bigg(I_0(\beta) - {\textbf L}_0(\beta) - 1 + \frac{2}{\pi}\,\mathrm{Si}(\beta)\bigg)\,.
\end{align*}
By numerical test it can be shown that for $\beta = m\,(\pi - \delta)$ it holds
$$
0 <  I_0(\beta) - {\textbf L}_0(\beta) - 1 + \frac{2}{\pi}\,\mathrm{Si}(\beta) < 1\,.
$$
Further it is known that $I_0(x) \geq 1$ for $x\in\R$, such that $0 < \Delta_1(\omega) < \frac{1}{I_0(\beta) - 1}$ for $\omega \in [- \delta,\,\delta]$.

Now we estimate $\Delta_2(\omega)$ for $\omega \in [-\delta,\,\delta]$ by the triangle inequality as
$$
	|\Delta_2(\omega)| \le \frac{\beta}{\pi\,\big(I_0(\beta) - 1\big)}\,\bigg(\int_1^{\nu_1(-\omega)} + \int_1^{\nu_1(\omega)}\bigg)\Big|\frac{\sin\big(\beta \sqrt{\nu^2-1}\big)}{\beta \sqrt{\nu^2-1}} - \frac{\sin(\beta \nu)}{\beta\,\nu}\Big|\,{\mathrm d}\nu\,.
$$
By \cite[Lemma~4]{PT21}, we have for $\nu \ge 1$ that
$$
	\Big|\frac{\sin\big(\beta\, \sqrt{\nu^2-1}\big)}{\beta\, \sqrt{\nu^2-1}} - \frac{\sin(\beta \nu)}{\beta \nu}\Big| \le \frac{2}{\nu^2}\,.
$$
Thus we conclude that
$$
	|\Delta_2(\nu)| \le \frac{4 \beta}{\pi\,\big(I_0(\beta) - 1\big)}\, \int_1^{\infty} \frac{1}{\nu^2}\,{\mathrm d}\nu = \frac{4 \beta}{\pi\,\big(I_0(\beta) - 1\big)}\,.
$$
Therefore we obtain for $\omega \in [-\delta,\,\delta]$ that
$$
	|\Delta(\omega)| \le \Delta_1(\omega) + |\Delta_2(\omega)| \le \frac{1}{I_0(\beta) -1}\,\bigg(1 + \frac{4  \beta}{\pi}\bigg)\,.
$$
Since the function ${\mathrm e}^{-x}\,I_0(x)$ is strictly decreasing on $[0,\, \infty)$ and tends to zero as $x\to \infty$ (see \cite{Ba10}), we have
$$
	\frac{1}{I_0(\beta) -1} = \frac{{\mathrm e}^{-\beta}}{{\mathrm e}^{-\beta}\,I_0(\beta) - {\mathrm e}^{-\beta}} = \frac{1}{{\mathrm e}^{-\beta}\,I_0(\beta) - {\mathrm e}^{-\beta}}\,{\mathrm e}^{-m\,(\pi - \delta)}\,.
$$
Thus, the approximation error of the continuous Kaiser--Bessel regularized Shannon formula~\eqref{eq:cKB} with the continuous Kaiser--Bessel window function \eqref{eq:RSS2} decreases exponentially with respect to~$m$.
\bend\\

Finally we show that for the continuous Kaiser--Bessel regularized Shannon sampling formula \eqref{eq:cKB} related to $\scr F_A$ and $f \in B_{p,\delta}^A$, the uniform perturbation error
$\| R_{\mathrm{cKB},m}^A f - R_{\mathrm{cKB},m}^A \tilde f \|_{\infty}$ only grows as ${\mathcal O}(m)$, if $\delta$ fulfills the oversampling condition $0 < \delta \leq \pi - \frac{\pi}{m}$.

\begin{theorem}\label{theo:cKBrobust}
For $m \in \N\setminus \{1\}$ and $\delta \in \big(0,\,\pi - \frac{\pi}{m}\big]$, let $f \in B_{p,\delta}^A$ be given.
Let $\varphi_{\mathrm{cKB}}\in\Phi_m$ be the continuous Kaiser--Bessel window function \eqref{eq:RSS2} with $\beta = m\,(\pi - \delta)$.
Furthermore, let $\tilde{f}(n)=f(n)+\varepsilon_n,\, n\in\Z$, be noisy samples,
where the complex sequence $\{\varepsilon_n\}_{n\in\Z}$ satisfies $|\varepsilon_n|\leq\varepsilon$ for some $\varepsilon>0$.\\
Then the continuous Kaiser--Bessel regularized Shannon sampling formula \eqref{eq:cKB} for $\scr F_A$ and $f \in B_{p,\delta}^A$ is numerically robust 
and 
$$
\|R^A_{\sinh,m}f - R_{\sinh,m}^A \tilde f\|_{\infty}\leq \varepsilon\, \big(2 + \frac{7}{4}\,  m\big)\,.
$$
\end{theorem}

\banf By Theorem \ref{theo:RSA2} we only have to estimate the value $\hat \varphi_{\mathrm{cKB}}(0)$ for the continuous Kaiser--Bessel window function \eqref{eq:RSS2}. From \eqref{eq:hatcKB} it follows that
\begin{eqnarray}
\hat \varphi_{\mathrm{cKB}}(0) &=& \frac{\sqrt 2\,m}{\sqrt \pi\,\big(I_0(\beta)-1\big)}\,\Big(\frac{\sinh \beta}{\beta} - 1 \Big)\\
&=& \frac{m\,{\mathrm e}^{\beta}}{\sqrt{2 \pi}\, \beta\,\big(I_0(\beta)-1\big)}\,\big(1 - {\mathrm e}^{-2 \beta} - 2 \beta\,{\mathrm e}^{-\beta}\big)\\
&<& \frac{m\,{\mathrm e}^{\beta}}{\sqrt{2 \pi}\, \beta\,\big(I_0(\beta)-1\big)}\,.
\end{eqnarray}
By \cite[p.~377, 9.7.1]{abst} one knows that
$$
\lim_{x \to \infty} \sqrt{2 \pi x}\,{\mathrm e}^{-x}\,I_0(x) = 1
$$
and hence also
$$
\lim_{x \to \infty} \sqrt{2 \pi x}\,{\mathrm e}^{-x}\,\big(I_0(x) - 1 \big) = 1\,.
$$
Numerical calculation shows that for all $x \geq \pi$ it holds
$$
\sqrt{2 \pi x}\,{\mathrm e}^{-x}\,\big(I_0(x) - 1 \big) \geq \sqrt 2 \, \pi\,\big(I_0(\pi) - 1 \big)\,.
$$
By the oversampling condition $\delta \in \big(0,\,\pi - \frac{\pi}{m}\big)$ we obtain $\beta = m\,(\pi - \delta) \geq \pi$.
Therefore we can estimate
$$
\frac{{\mathrm e}^{\beta}}{\sqrt{2 \beta}\,\big(I_0(\beta)-1\big)} \leq \frac{{\mathrm e}^{\pi}}{\sqrt{2 \pi}\,\big(I_0(\pi)-1\big)} = 1.163167956\,.
$$
Since $1.163167956 \times \sqrt 2 < \frac{7}{4}$, the above estimate of the perturbation error is shown. \bend

\section{Numerical Comparison of Proposed Regularization Methods}

In this section, we illustrate how the regularization techniques proposed in this paper perform in numerical practice when reconstructing SAFT-bandlimited functions both in the noise-free and in noisy samples scenario.  %from the stability perspective. 
We start with describing a test function of your choice for the numerical experiments and then move on to the presentation of numerical results.

\begin{figure}[ht]
    \centering 
    \subfloat[\label{fig:-1}]{
        \includegraphics[width=0.45\linewidth]{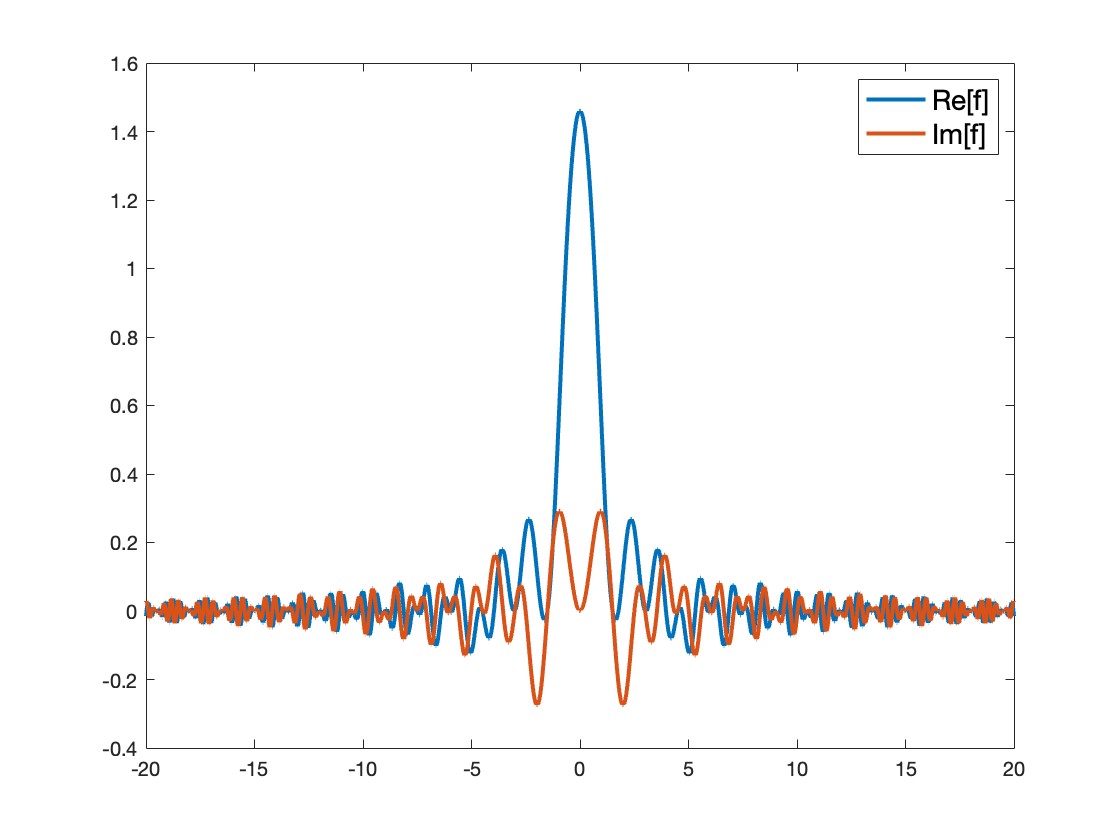}
    }
    \subfloat[\label{fig:-2}]{
       \includegraphics[width=0.45\linewidth]{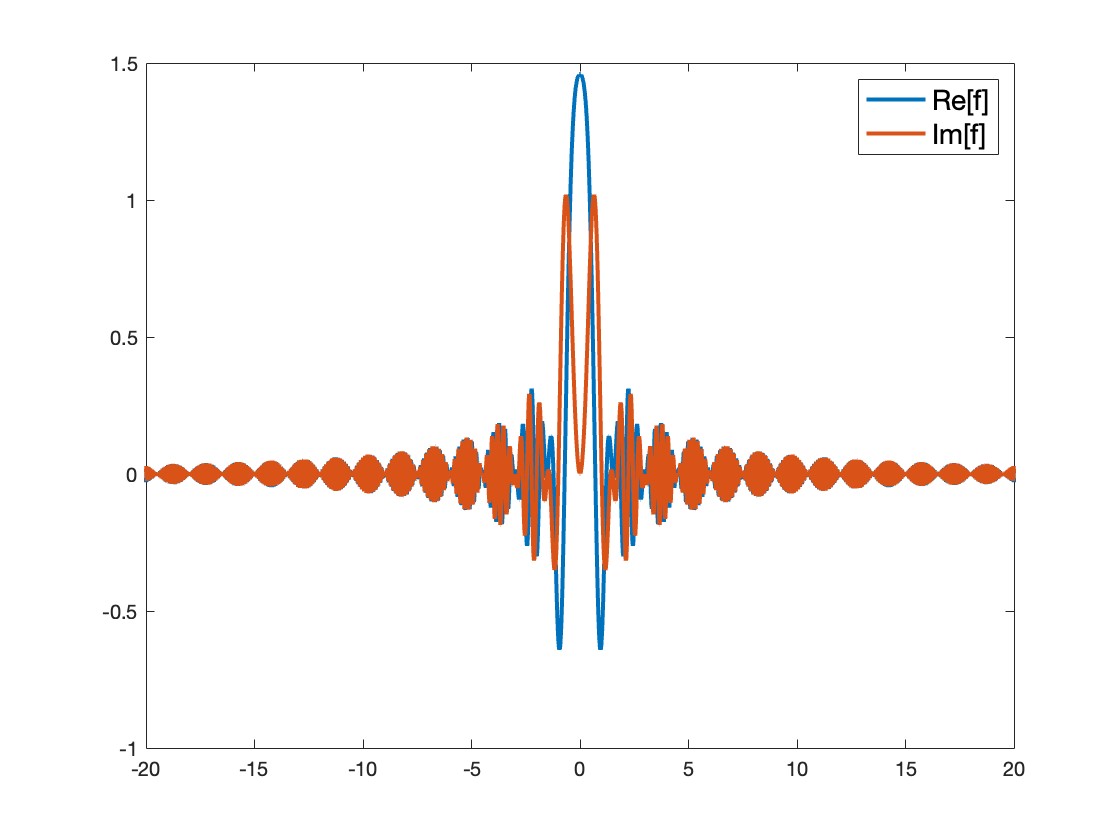}
    }
    \caption[]{ The test function $f$ with fixed parameter $h=1.5,$ $ a=cos(\alpha), b=sin(\alpha)$ for $\alpha=\frac{\pi}{4}$ (right)  and  $\alpha=\frac{\pi}{20}$ (left). The smaller the angle $\alpha$ is, the the stronger the function oscillates. }
    \label{fig:1:num}
\end{figure}

To construct a function in $f\in B_{p,\delta}^A$, we proceed as follows. We start with a bandlimited function $g\in B_{\delta}$ with $\delta=\pi/h\le \pi$ and $h\in (1,2]$, here the condition $\delta<\pi$ ensures the oversampling process
 for $f$. For parameter $k\in\N $ and $h\in \R$, we consider a family of the functions $S_{k,h}$ defined as
 \begin{equation*}
     S_{k,h}(t)=\mathrm{sinc}\big((t-k h)/h\big), \quad k\in \Z.
 \end{equation*}
According to \cite{S93}, for each $h$, functions $S_{k,h}$ form an orthogonal system in $B_{\delta}$, with respect for $k$, namely it holds that 
  \begin{equation*}
    \int_{\R} S_{k,h}(t)S_{\ell,h}(t)\dx t= h \delta_{k,\ell}, \quad k, \ell \in \Z.
 \end{equation*}
We choose as an auxiliary function, the function $g$ defined as the following sum 
 \begin{equation*}
     g(t)=\frac{2}{\sqrt{5h}}\big(S_{0,h}(t)+ \frac{1}{2}S_{1,h}(t)\big), \quad t\in \R.
 \end{equation*}
 It is easy to check that $g$ is normalized, that is $\int_{\R}g^2(t)\dx t=1$. Using $g$, the test function $f$ of our choice will be defined by the following relation
 \begin{equation}\label{eq:test-fun}
     f(t)=\e^{\im\frac{a}{2b}t^2} g(t),
 \end{equation}
which according to Proposition~\ref{pro:SSS1} belongs to the space $B_{p,\delta}^A$ and has the unit $L^2(\R)$-norm.
The real and imaginary parts of $f$ are depicted in Figure~\ref{fig:1:num} for different parameter choices.

For the experimental setting, we fix the interval length to $N=50$, and consider the reconstruction of  the test function, on the interval $[-N,N]$
from $2N+2m-1$ samples $f(n)$ with $n=1-N-m,\,\dots,\, N+m-1$. Having this data, we compare the performance of the classical reconstruction technique \eqref{eq:SSS1} to the reconstruction technique \eqref{eq:RSS0} proposed here, exploring numerically how the reconstruction error behaves when the time window size $m$ grows.  Herewith, we consider the regularization with three different window functions,  the centered cardinal B-spline, sinh-type, and Kaiser-Bessel window, all indicated in~Example~\ref{ex:RSS1}. 

In the first part of our experiments,  we consider the noise-free scenario and measure the maximum approximation error 
\begin{equation}\label{eq:error_no_noise_num_sec}
\max_{t \in [-N,\,N]} \big|f(t) - (R_{\varphi,m}^A f)(t)\big|\,,
\end{equation}
where the reconstructed version $R_{\varphi,m}^A f$ of the function $f\in B_{p,\delta}^A$ is defined as
\begin{equation}\label{eq:recon_of_f_num_sec}
(R_{\varphi,m}^A f)(t) = \sum_{n=1-N-m}^{N+m-1} f(n)\,{\mathrm e}^{{\mathrm i}\frac{a}{2b}\,(t^2 - n^2)}\,\mathrm{sinc}(t-n)\,\varphi(t-n)\,, \quad t \in [-N,\,N]\,, 
\end{equation}
with some window function  $\varphi \in \Phi_m$. Following the results in Theorem~\ref{theo:SRSS1}, \ref{theo:SRRS2} and \ref{theo:cKB}, we consider the three corresponding cases with $\varphi \in \{\varphi_{\mathrm B},\, \varphi_{\sinh},\,\varphi_{\mathrm{cBK}}\}$, and in our experiments we let the window size $ m$ to settlement in the set $ m \in \{11+3n, n=1,2,3,5,6 \}$. To approximate \eqref{eq:error_no_noise_num_sec}, we evaluate the difference $f - (R_{\varphi,m}^A f)$ at $10^5$ equispaced point in the interval $[-50, 50]$.

\begin{figure}[ht!]
    \centering 
    \subfloat[\label{fig:-1}]{
        \includegraphics[width=0.45\linewidth]{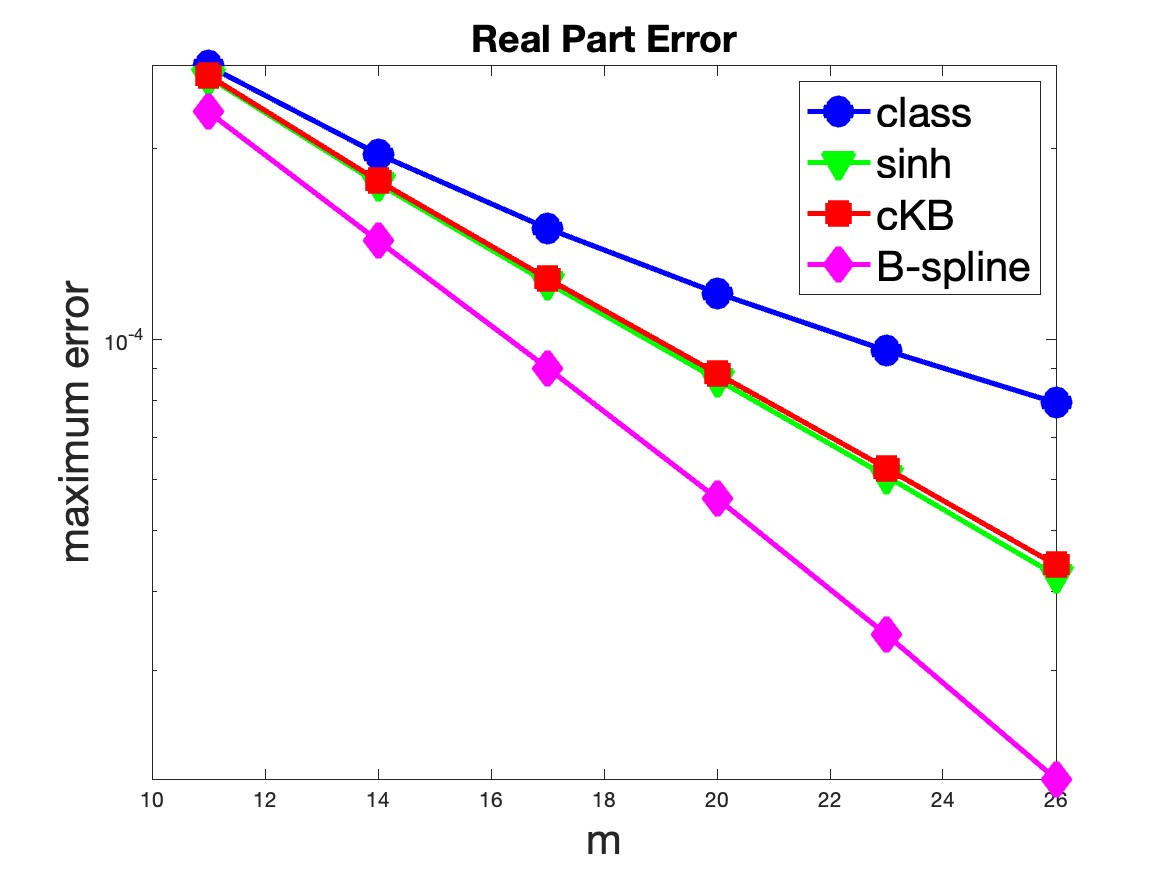}
    }
    \subfloat[\label{fig:-2}]{
       \includegraphics[width=0.45\linewidth]{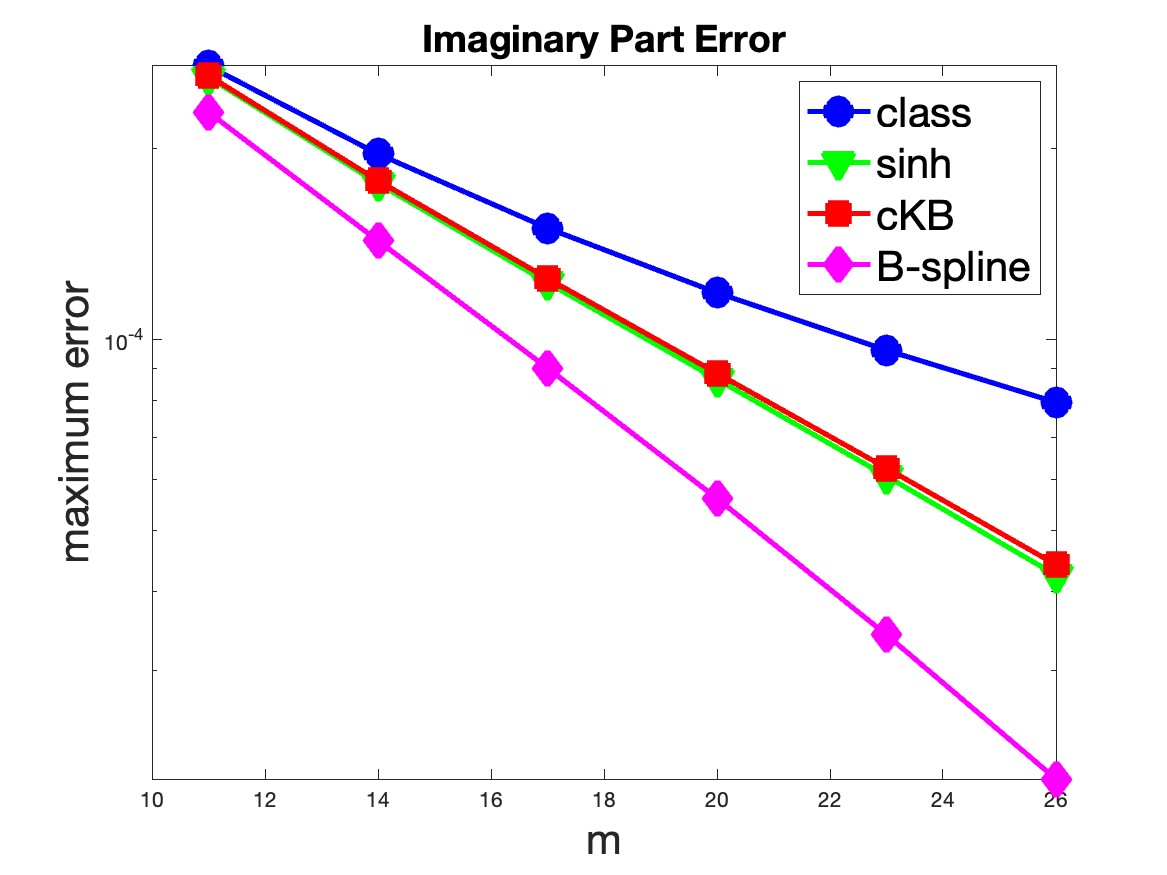}
    }
    \caption[]{Maximum approximation error on the interval $[-50,50]$ by the reconstruction of the test function~\eqref{eq:test-fun} from its noiseless samples via the classical reconstruction technique \eqref{eq:SSS1} and the regularized one \eqref{eq:RSS0} using time window functions $\varphi_{sinh}$, $\varphi_{cKB}$ and $\varphi_{B}$ with the window size $ m \in \{11+3n, n=1,2,3,5,6 \}$ }
    \label{fig:1:num}
\end{figure}

The results of this numerical experiment are depicted in Figure~\ref{fig:1:num}. As one can observe, for all the methods the maximum error decays when the window size $m$ and so the number of samples increases. All the regularization techniques proposed here outperform the classical reconstruction technique and the performance is obtained with the centered cardinal B-spline window $\varphi_{B}$. At the beginning, the error for B-splines is of order $2.28\cdot 10^{-4}$ and then it decays reaching its lowest point at $2.023\cdot 10^{-5}$.

The second part of the experiments considers the case when the samples of $f$ are corrupted by noise. In this case, instead of $f(n)$, the obtained samples are given as $\tilde{f}(n)=f(n)+\varepsilon_n,\, n\in\Z$
with some complex-valued noise term  $\varepsilon_n$, for $2N+2m-1$. In our experiment, for each $n$, we take the elements $\varepsilon_n$ as independent complex random variables with real and imaginary parts, $\Re({\varepsilon_n})$ and $\Im({\varepsilon_n})$, continuous uniformly distributed in the interval $10^{-5}\cdot[1,\, 5]$.
From the results in Theorem~\ref{theo:Brobust},  \ref{theo:sinhrobust} and 
\ref{theo:cKBrobust}, we consider 
\begin{eqnarray}\label{eq:error_with_noise_num_sec}
\max_{t \in [-N,\,N]} \big|f(t) - R^A_{\varphi,m}\tilde{f}(t)\big|\,,&
\end{eqnarray}
where  for $(R_{\varphi,m}^A \tilde{f})(t) $ the samples ${f}(n)$ in the formula \eqref{eq:recon_of_f_num_sec}  are replace by its noisy counterparts $\tilde{f}(n)$. Similarly, to the noise-free scenario, we approximate \eqref{eq:error_with_noise_num_sec} by evaluating the difference $f - (R_{\varphi,m}^A \tilde f)$ at $10^5$ equispaced point in the interval $[-50, 50]$. Moreover, here for each $ m \in \{11+3n, n=1,2,3,5,6 \}$, we run $100$ experiments, and as the final performance, we consider the average of maximum error.

\begin{figure}[ht]
    \centering 
    \subfloat[\label{fig:-1}]{
        \includegraphics[width=0.45\linewidth]{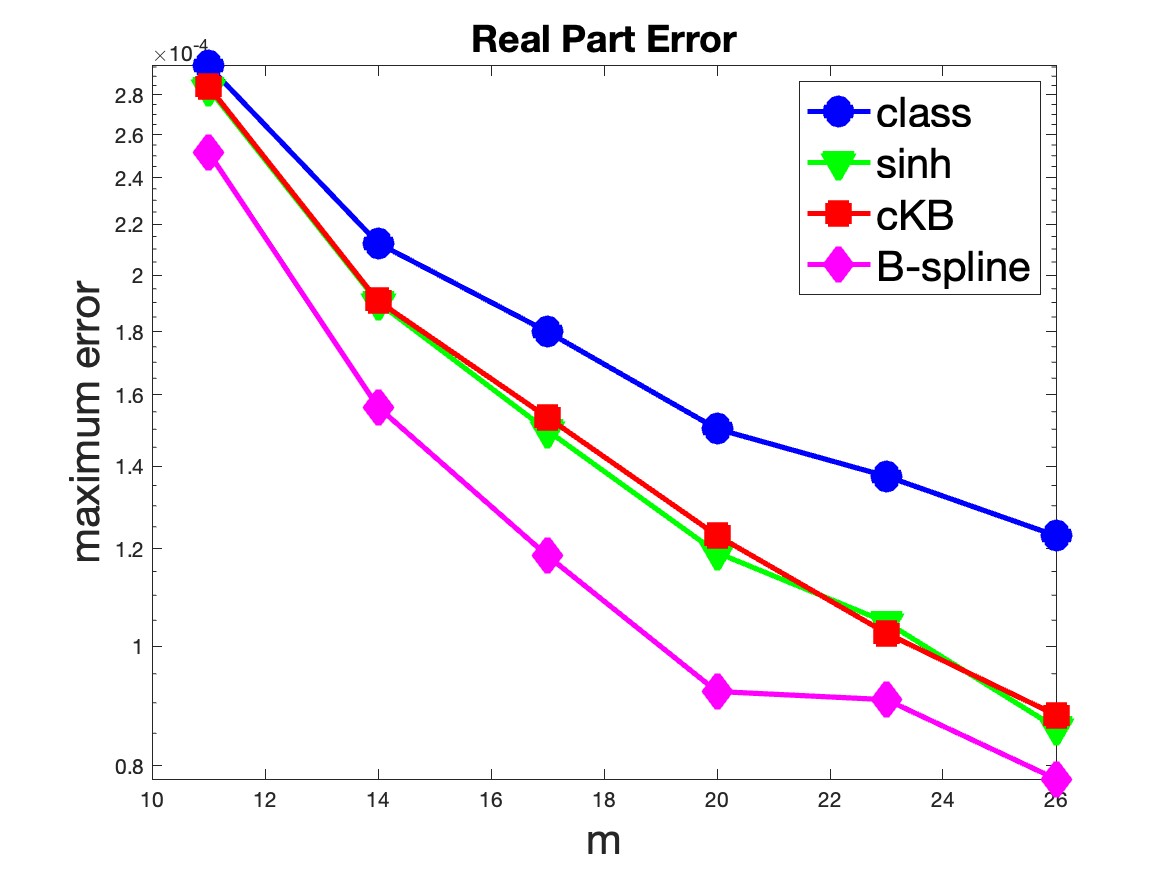}
    }
    \subfloat[\label{fig:-2}]{
       \includegraphics[width=0.45\linewidth]{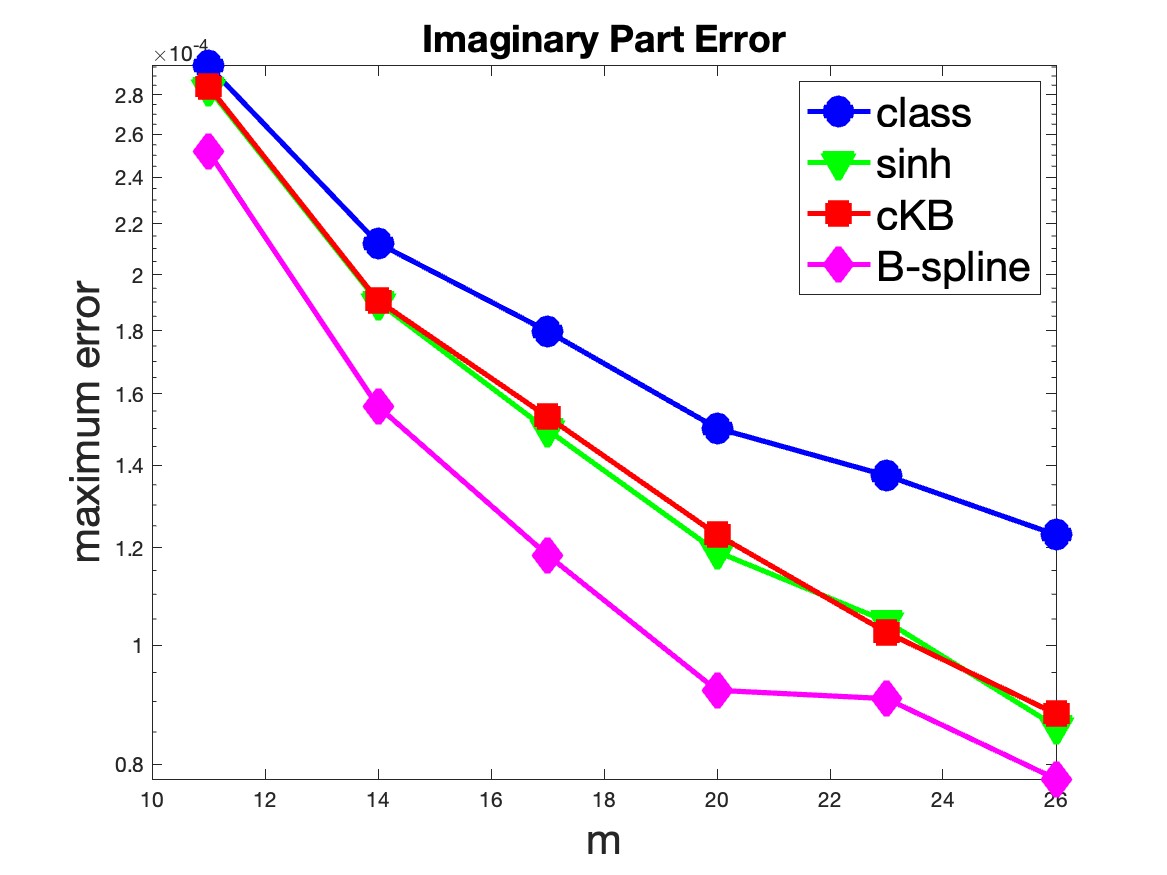}
    }
    \caption[]{Maximum approximation error on the interval $[-50,50]$, averaged over 100 runs, by the reconstruction of the test function~\eqref{eq:test-fun} from the samples $f(n)+\varepsilon_n$ corrupted by complex-valued noise $\varepsilon_n$ with $|\varepsilon_n|\leq10^{-4}$, $n=1-N-m,\,\dots,\, N+m-1$,  via the classical reconstruction technique \eqref{eq:SSS1} and the regularized one \eqref{eq:RSS0} using time window functions $\varphi_{sinh}$, $\varphi_{cKB}$ and $\varphi_{B}$ with the window size~${ m \in \{11+3n, n=1,2,3,5,6 \}}$ }
    \label{fig:2:num}
\end{figure}

% $$
% (R_{\varphi,m}^A \tilde{f})(t) = \sum_{n=1-N-m}^{N+m-1} \tilde{f}(n)\,{\mathrm e}^{{\mathrm i}\frac{a}{2b}\,(t^2 - n^2)}\,\mathrm{sinc}(t-n)\,\varphi(t-n)\,, \quad t \in [-N,\,N]\,,
% $$
Figure~\ref{fig:2:num} shows the results of numerical experiments with noisy data. The proposed regularization methods perform better in the noisy case too as similar to the noise-free case, the maximum error decays when the window size $m$ increases. Also here, the centered cardinal B-spline window $\varphi_{B}$ allows for the smallest maximum error. At the beginning, the  B-spline error is of order $2.5\cdot 10^{-4}$ and decays to $7.78\cdot 10^{-5}$, remaining within the order of perturbation.

In conclusion, we can observe that the regularized Shannon sampling formulas for SAFT-bandlimited functions exhibit superior performance over the classical Shannon sampling formulas not only in the noise-free case but also when the data are corrupted by noise. 

\newpage

%
%
%<<<<<<<<<<<<<<<<<<<<<<<<<<<<<<<<<<<<<< References >>>>>>>>>>>>>>>>>>>>>>>>>>>>>>>>>>>>>>>>>>>>>>>>>>>>>>>>>>
%

%
%<<<<<<<<<<<<<<<<<<<<<<<<<<<<<<<<<<<<<<<<<<<<<<<<<>>>>>>>>>>>>>>>>>>>>>>>>>>>>>>>>>>>>>
%
\end{document}